\input amstex

\documentstyle{amsppt}
\loadbold
\loadmsbm
\UseAMSsymbols

 \NoRunningHeads
 \NoBlackBoxes

\topmatter
\title
Irreducible Koopman representations for nonsingular  actions on~boundaries of rooted trees
\endtitle

\author  Alexandre I. Danilenko and Artem Dudko
\endauthor

\address
Faculty of Mathematics and Computer Science,
Nicolaus Copernicus University,
Chopin street 12/18, 87-100 Toru{\'n}, Poland
\newline\indent
 B. Verkin Institute for Low Temperature Physics and Engineering
of the  National Academy of Sciences of Ukraine,
47 Nauky Ave.,
 Kharkiv, 61164, Ukraine
  \newline
 \indent
 Mathematical Institute of the Polish Academy of Sciences,
ul. {\'S}niadeckich 8,
00-656,
 Warszawa,
 Poland
\endaddress
\email            alexandre.danilenko\@gmail.com
\endemail

\address
Mathematical Institute of the Polish Academy of Sciences,
ul. {\'S}niadeckich 8,
00-656,
 Warszawa,
 Poland
\endaddress
\email    adudko\@impan.pl
\endemail

\thanks
This work  was  supported in part by the
``Long-term program of
support of the Ukrainian research teams at the
Polish Academy of Sciences carried out in
collaboration with the U.S. National Academy of
Sciences with the financial support of external
partners''.
\endthanks

\abstract
Let $G$ be a countable branch group of automorphisms of a spherically homogeneous rooted tree.
Under some  assumption on finitarity of $G$,
we construct,
for each sequence $\omega\in\{0,1\}^\Bbb N$,  an irreducible unitary representation
$\kappa_\omega$ of $G$.
Every two representations $\kappa_\omega$ and $\kappa_{\omega'}$ are
weakly equivalent.
They are unitarily equivalent if and only if
$\omega$ and $\omega'$ are tail equivalent.
Each $\kappa_\omega$ appears as the Koopman representation associated with some ergodic $G$-quasiinvariant measure (of infinite product type) on the boundary of the tree.
\endabstract

\endtopmatter

\document

\head 0. Introduction
\endhead

A spherically homogeneous rooted tree $T$
 is a graded graph whose
 vertices  are partitioned into levels numbered with non-negative integers.
 The graph is parametrized by a sequence $(q_n)_{n=1}^\infty$ of integers $q_n\ge 2$.
There is a single vertex on the $0$-th level referred to  as the root of $T$.
There are $q_1\cdots q_n$ vertices on the $n$-th level.
Each vertex of the $n$-th level is connected with a single vertex on the $(n-1)$-th level and with
$q_{n+1}$ vertices on the $(n+1)$-th level.
Denote by $\partial T$ the boundary of $T$.
The boundary consists of all infinite simple paths in $T$ starting at the root (i.e. $0$-th level).
We identify $\partial T$ naturally with the infinite product $\prod_{n=1}^\infty\{1,\dots, q_n\}$.
We endow $\partial T$ with the natural product of discrete topologies on the factors
$\{1,\dots, q_n\}$.
Every automorphism of $T$ (we assume that the automorphisms preserve the levels of $T$)
induces a continuous transformation of $\partial T$.
Hence, each countable group $G$ of automorphisms of $T$ induces a continuous $G$-action on $\partial T$.
In this paper we consider only those $G$  whose associated actions on $\partial T$ are minimal.
Our purpose is to construct  a large  family of mutually  inequivalent irreducible
unitary representations of $G$.
This will be done for the class of branch groups (see Definition~1.7) satisfying certain finitarity condition (see \S3 for definition of $\Cal A_\lambda^\bullet$).

Branch groups include many interesting examples with unusual properties.
One of the first and most important examples of branch groups was introduced by R.~Grigorchuk in \cite{Gr1} and \cite{Gr4} as finitely generated torsion group of intermediate growth.
Examples of branch groups appear naturally as iterated monodromy groups in dynamics of rational maps \cite{Ne}.
One of them is the iterated monodromy group of the map $z\mapsto z^2+i$ (see \cite{BaGrNe}).
 Also, branch groups play important role in classification of infinite groups with respect
  to the property of {\it just-infiniteness} \cite{Gr5}.
%. It was shown to be amenable but not subexponentially amenable (see \cite{BaVi} and  \cite{GrZh}).

Koopman representations of groups acting on rooted trees, including branch and weakly branch groups, were studied, for instance, in \cite{BaGr} and \cite{DuGr1}. There, the authors obtained various results on questions of ergodicity of the actions of these groups, irreducibility and disjointness of Koopman representations, their weak equivalence, and other topics. In this paper we present a different approach to these problems. This approach relies on techniques from non-singular dynamics. Therefore, we use the terminology and notations which are common to the field of non-singular dynamics rather then to the area of groups acting on rooted trees. We hope that this will be useful for both the readers interested in actions on rooted trees and non-singular actions.

%{\color{blue} The main result concerns the class of groups acting \emph{purely $\lambda$-finitarily} on a spherically homogeneous rooted tree $T$, where $\lambda$ is the uniform Bernoulli measure $\}
\proclaim{Main Result  (see Theorem~6.6)} Let $G$ be a branch subgroup of $\Cal A_\lambda^\bullet$.
Then there is a continuous mapping $\omega\mapsto\mu_\omega$ from the set
$\{0,1\}^\Bbb N$ to the set of infinite product probabilities on the boundary $\partial T$ of the spherically homogeneous rooted tree $T$
such that:
\roster
\item"(i)" $G$ is nonsingular and ergodic with respect to $\mu_\omega$ for each $\omega\in \{0,1\}^\Bbb N$.
\item"(ii)" The Koopman representation $\kappa_{\mu_\omega}$ of $G$ associated with the dynamical system
$(\partial T,\mu_\omega,G)$ is irreducible for each $\omega\in \{0,1\}^\Bbb N$.
\item"(iii)"
If $\omega,\omega'$ from $\{0,1\}^\Bbb N$ are tail equivalent then
$\kappa_{\mu_\omega}$ and $\kappa_{\mu_\omega}$ are unitarily equivalent.
\item"(iv)"
If $\omega,\omega'$ from $\{0,1\}^\Bbb N$ are not tail equivalent then
$\kappa_{\mu_\omega}$ and $\kappa_{\mu_\omega}$ are disjoint (not unitarily equivalent).
\item"(v)"
$\kappa_{\mu_\omega}$ and $\kappa_{\mu_\omega}$ are weakly equivalent
for all $\omega,\omega'\in \{0,1\}^\Bbb N$.
\endroster
\endproclaim

Thus, we obtain an uncountable family of mutually non-unitarily equivalent classes of irreducible unitary representations of $G$.
We may consider this family ``complicated'' or ``unclassifiable''  because there exists no  Borel subset $B$ of
$ \{0,1\}^\Bbb N$  such that the set $\{\kappa_{\mu_\omega}\mid \omega\in B\}$
meets each unitary equivalence class exactly once.

The problem of finding  uncountable families of irreducible Koopman representations for groups of $T$-automorphisms
 was considered earlier in \cite{DuGr1} in a particular case of a regular rooted tree (when $q_1=q_2=\cdots$) and $\partial T$ endowed with Bernoulli measures $\beta^{\otimes\Bbb N}$, where $\beta$ is a probability on $\{1,\dots,q_1\}$.
 The authors of that paper looked for conditions for $G$ that guarantee  ergodicity of $G$ with respect to any element of the entire simplex of Bernoulli measures on $\partial T$.
 It was stated in \cite{DuGr1, Proposition~3} that if $G$ is subexponentially bounded and minimal
 then it is ergodic for each Bernoulli measure.
Unfortunately, there was a flaw in the proof of the above claim.
 We construct counterexamples to \cite{DuGr1, Proposition~3} in \S4 below.
Thus, the statement of irreducibility of representations in \cite{DuGr1, Theorem~3,1)} and
\cite{DuGr2, Corollary~23} which is based  on \cite{DuGr1 Proposition~3}, requires additional conditions to ensure ergodicity of the group action. 
We present an example of such conditions in \S5.

In our opinion,  there exist no ``universal'' conditions on $G$ that imply ergodicity
for the entire class of Bernoulli measures.
However,
since ergodicity is necessary  for  irreducibility of the Koopman representations,
we develop another approach to this problem.
Instead of  fixing a class of measures (Bernoulli measures) and looking for  transformation groups ``compatible'' with those measures, we start with a group $G$ of  automorphisms of the rooted tree $T$ and then select a class of measures ``compatible'' with this group.
If $G$
is branch and satisfies  certain  finitarity condition, we can find a suitable (depending on $G$) uncountable  collection of compatible measures within the family
of
infinite product measures on $\partial T$.

Our proof of unitary equivalence and disjointness of the Koopman representations is different from
\cite{DoGr1}.
It is based on the Kakutani dichotomy theorem on equivalence and orthogonality of infinite product measures \cite{Ka}.
%``Most'' of the factors of these measures, i.e. the pushforwards  onto $\{1,\dots,q_n\}$,
%are the equidistributions.
%The ``remaining''  factors are chosen  so that the infinite product of all factors is
% non-atomic, $G$-qusiinvariant and ergodic.

The outline of the paper is as follows.
In Section~1 we consider infinite wreath products of finite groups and discuss their relations with groups of automorphisms of spherically homogeneous rooted trees.
We isolate minimal subgroups of such automorphisms, list their basic properties and show that these groups can be interpreted as odometer actions.
In Section~2 we remind  some basic concepts of nonsingular dynamics and show in Proposition~2.1 that
the weakly branch groups are ``universally''  (i.e. with respect to every quasiinvariant measure) rigid.
In Section~3 we study various ``finitarity'' conditions on $G$ that imply nonsingularity of the action of $G$
 with respect to a given infinite product measure on $\partial T$.
 Examples~3.6 and 3.7 illustrate a subtle dependence of these conditions on the underlying measure.
Section~4 is devoted to the concept of ``subexponential boundedness'' introduced in \cite{DuGr1}.
 In Example~4.3, using the Shannon-McMillan theorem, we show that our concept of $\mu$-finitarity is more general than the subexponential boundedness.
 %This answers a question by R.~Grigorchuk.
We also produce 3 examples of subexponentially bounded minimal groups $G$ which are not ergodic with respect to Bernoulli measures on $\partial T$:
$G$ is dissipative in Example~4.4, $G$ is conservative non-ergodic in Example~4.5 and $G$ is conservative non-ergodic weakly branch in Example~4.6.
A way to construct uncountable collections of  infinite product measures on $\partial T$ which are compatible with finitary groups of $T$-automorphisms
is elaborated in Section~5: see Theorem~5.3 and Corollary~5.4.
In Section~6 we consider Koopman representations associated with nonsingular $G$-actions and compatible measures.
%The concept of quasi-branch group is introduced there.
The main result of the paper (except for (v)) is proved in Theorems~6.4 and~6.5.
The final Section~7 is devoted to the weak equivalence of Koopman representations.
We prove there part (v) of the main result.

{\it Acknowledgements.}  The authors are grateful to Rostislav Grigorchuk for his 
valuable comments.

\head 1. Minimal  actions on the boundaries of spherically homogeneous rooted trees
\endhead

\subhead 1.1. Infinite  wreath products
\endsubhead
Let $A$ and $H$ be finite groups and let $H$ act on a finite set $\Omega$.
We recall that  the regular wreath product $A\wr_\Omega H$ is the semidirect product $A^\Omega\rtimes H$,
where $H$ acts on $A^\Omega$ in a natural way, i.e.
$$
h\cdot (a_\omega)_{\omega\in\Omega}:=(a_{h^{-1}\omega})_{\omega\in\Omega}.
$$
Hence, a short exact sequence
$$
1\to A^\Omega\to A\wr_\Omega H\to H\to 1
$$
is well defined.
We note that the projection $p: A\wr_\Omega H\to H$ in this sequence splits, i.e. there is a homomorphism
$i:H\to A\wr _\Omega H$ such that $p\circ i=\text{id}$.

Suppose now that we have an infinite  sequence $(A_n)_{n=1}^\infty$ of finite groups and $A_n$ acts
on a certain finite space $\Omega_n$ for each $n\in\Bbb N$.
Then $A_2\wr_{\Omega_1}A_1$ acts naturally on $\Omega_2\times \Omega_1$.
Hence, the group $A_3\wr_{\Omega_2\times\Omega_1}(A_2\wr_{\Omega_1}A_1)$ is well defined.
This group, in turn, acts naturally on $\Omega_3\times \Omega_2\times \Omega_1$, and so on.
By the {\it infinite regular wreath product\,}
$$
\cdots\wr_{\Omega_3\times \Omega_2\times \Omega_1} A_3 \wr_{\Omega_2\times \Omega_1} A_2\wr_{\Omega_1} A_1
$$
 we mean the projective limit of the sequence
$$
 \cdots\to A_4\wr_{\Omega_3\times\Omega_2\times\Omega_1}(A_3\wr_{\Omega_2\times\Omega_1} (A_2\wr_{\Omega_1} A_1))
 \to A_3\wr_{\Omega_2\times\Omega_1} (A_2\wr_{\Omega_1} A_1)\to     A_2\wr_{\Omega_1} A_1\to A_1.
%
 %\leftarrow  A_4\wr(A_3\wr (A_2\wr A_1))\leftarrow\cdots.
$$
On the other hand, as every arrow in this sequence splits, we have  a sequence of group embeddings:
$$
A_1@>{i_1}>> A_2\wr_{\Omega_1} A_1@>{i_2}>> A_3\wr_{\Omega_2\times\Omega_1} (A_2\wr A_1)@>{i_3}>>\cdots.
$$
We call the inductive limit of this sequence {\it the infinite restricted wreath product of $(A_n)_{n=1}^\infty$}.
Denote it by $\Cal A_0$.
It is a countable group.
As the diagram
$$
\CD
 %\cdots\to A_4\wr(A_3\wr (A_2\wr A_1))
  \cdots @>>> A_3\wr_{\Omega_2\times\Omega_1} (A_2\wr_{\Omega_1} A_1)@>>> A_2\wr_{\Omega_1} A_1@>>> A_1\\
  \cdots @.@AA{i_2}A @AA{i_1}A @AAA\\
    \cdots @<{i_2}<< A_2\wr_{\Omega_1} A_1@<{i_1}<< A_1@<<<1
 % \cdots\to A_4\wr(A_3\wr (A_2\wr A_1))\to A_3\wr (A_2\wr A_1)\to A_2\wr A_1\to A_1.
\endCD
$$
commutes, it defines
 an embedding of $\Cal A_0$ into $\cdots\wr_{\Omega_3\times\Omega_2\times\Omega_1} A_3 \wr_{\Omega_2\times\Omega_1}A_2\wr_{\Omega_1} A_1$.

Endow the infinite regular wreath product with the natural infinite product of discrete topologies.
Then $\cdots\wr_{\Omega_3\times\Omega_2\times\Omega_1} A_3 \wr_{\Omega_2\times\Omega_1}A_2\wr_{\Omega_1} A_1$ is 0-dimensional compact topological group.
The infinite restricted wreath product of $(A_n)_{n=1}^\infty$ is a dense countable subgroup in
the infinite regular wreath product.
Thus, we consider $\cdots\wr_{\Omega_3\times\Omega_2\times\Omega_1} A_3 \wr_{\Omega_2\times\Omega_1}A_2\wr_{\Omega_1} A_1$ as the natural 0-dimensional group compactification of $\Cal A_0$.

\subhead 1.2. Dynamical realization of  infinite wreath products
\endsubhead
Fix a sequence $(q_n)_{n=1}^\infty$ of natural numbers.
We will assume that $q_n\ge 2$ for each $n$.
Let $X_n:=\{1,\dots,q_n\}$ and let
$X:=X_1\times X_2\times\cdots$.
Endow $X_n$ with the discrete topology for each $n$ and endow $X$ with the standard infinite product topology.
Then $X$ is a compact Cantor set.
For $m>n>0$, we let
 $X_n^m:=X_n\times\cdots\times X_m$.
 Given $(x_1,\dots,x_n)\in X_1^n$, we denote by $[x_1,\dots,x_n]_1^n$ the corresponding $n$-cylinder in $X$:
$$
[x_1,\dots,x_n]_1^n:=\{(y_j)_{j=1}^\infty \in X\mid y_j=x_j, 1\le j\le n\}.
$$
Every cylinder is a clopen subset of $X$.
Every clopen subset of $X$ is a union of finitely many mutually disjoint cylinders.
Denote by $\pi_n$ the projection mapping $\pi_n:X\ni (x_j)_{j=1}^\infty\mapsto(x_j)_{j=1}^n\in X_1^n$.

\definition{Definition 1.1}
We say that a homeomorphism $g:X\to X$ is {\it compatible with} $(\pi_n)_{n=1}^\infty$ if
 for each $n>0$, there is a bijection $g_n$ of $X_1^n$ such that $\pi_n\circ g=g_n\circ\pi_n$.
\enddefinition

The following proposition is straightforward.
We leave its proof to the reader.

\proclaim{Proposition 1.2} Let $g$ be a homeomorphism of $X$.
Then $g$ is compatible with $(\pi_n)_{n=1}^\infty$ if and only if
there exist a bijection $\gamma$ of $X_1$ and a sequence of mappings $\alpha_n:X_1^n\to\text{\rom{Bij}}(X_{n+1})$, $n\in\Bbb N$ such that
$$
g(x_1,x_2,\dots):=(\gamma[x_1], \alpha_1(x_1)[x_2], \alpha_2(x_1,x_2)[x_3],\dots).
\tag 1-1
$$
\endproclaim

Denote by $\Cal A$ the infinite regular wreath product
$$
\cdots\wr_{X_3\times X_2\times X_1}\text{Bij}(X_3)\wr_{X_2\times X_1}\text{Bij}(X_2)\wr_{X_1}\text{Bij}(X_1).
$$
The infinite restricted wreath product of $(\text{Bij}(X_n))_{n=1}^\infty$ will be denoted by
$\Cal A_0$.
It follows from Proposition~1.2  that the group of homeomorphisms of $X$ which are compatible with
$(\pi_n)_{n=1}^\infty$
is isomorphic to $\Cal A$.
The formula \thetag{1-1} determines an embedding of $\Cal A$ into Homeo$(X)$.
We interpret \thetag{1-1} as a dynamical realization of $\Cal A$ by homeomorphisms of a Cantor set.

\remark{Remark 1.3} Let  $T$ be a   spherically homogeneous rooted tree.
Then  $\partial T$ is identified naturally with $X$.
The group of automorphisms of $T$  (that preserve the levels) is exactly $\Cal A$.
This group acts naturally on $\partial T$.
Of course, this action is given by~\thetag{1-1}.
\endremark

\subhead 1.3.
Minimal subgroups of $\Cal A$
\endsubhead
From now on we consider $\Cal A$ as a group of transformations of $X$.
Given a subgroup $G\subset\Cal A$, we let $G_n:=\{\pi_n\circ g\mid g\in G\}$.
Then $G_n$ is a  subgroup of Bij$(X_1^n)$.
The following claim is straightforward.

\proclaim{Claim 1.4 \cite{GrNeSu, Proposition 6.5}}  The following are equivalent:
\roster
\item"(i)" The $G$-action on $X$ is minimal.
\item"(ii)"
The $G$-action on $X$ is topologically transitive.
\item"(iii)"
$G_n$ is transitive on $X_1^n$ for each $n\in\Bbb N$.
\endroster
\endproclaim

\proclaim{Claim 1.5 \cite{GrNeSu, Proposition 6.5}}  The action of $G$ on $\partial T$ is minimal if and only if it is uniquely ergodic.
\endproclaim
\demo{Proof}
Let $G$ be minimal
and  let $\mu$ be a $G$-invariant probability on $X$.
We set $\mu_n:=\mu\circ\pi_n^{-1}$.
Then $\mu_n$ is a $G_n$-invariant probability on $X_1^n$.
Since $G_n$ is transitive by Claim~1.4, it follows that
$\mu_n$ is the equidistribution on $X_1^n$.
Since $\mu:=\projlim_{n\to\infty}\mu_n$,
we obtain that $\mu$ is defined uniquely.
Thus, $G$ is uniquely ergodic.

Conversely, let $G$ be uniquely ergodic.
Denote by $\lambda$ the projective limit of the equidistributions $\lambda_n$ on $X_n$.
Then $\lambda$ is a $G$-invariant probability on $X$.
Hence $\lambda$ is ergodic.
Hence $G_n$ is ergodic on $(X_1^n,\lambda_n)$ for each $n$.
This is only possible if $G_n$ is transitive on $X_1^n$.
It remains to apply Claim~1.4.
\qed
\enddemo

We note that if $G$ is minimal then the unique $G$-invariant probability measure on $X$ does not depend on the choice $G$.
We call this measure {\it the Haar measure for $X$}.
It is the infinite direct product of the equidistributions on $X_n$, $n>0$.

Consider now the case of qusiinvariant measures.
The following is a partial analogue of Claim~1.5.

\proclaim{Claim 1.6} Let $G$ be a subgroup in $\Cal A$ and let $\mu$ be a $G$-quasiinvariant finite measure on $X$.
If $G$ is $\mu$-ergodic
and  and $\mu(O)>0$ for each non-empty clopen subset $O$ of $X$ then $G$ is minimal.
%\item"(ii)"
%If $G$ is minimal and the tail $\sigma$-algebra on $X$ is trivial (mod $\mu$) then $G$ is $\mu$-ergodic.
%\endroster
\endproclaim

It should be noted that the converse to Claim~1.6 does not hold: there exist minimal transformation group $G$ and a non-ergodic $G$-quasiinvariant measure with $\mu(O)>0$ for each non-empty clopen subset $O$ of $X$.
See Examples~3.6--3.8 below.

Taking into account Remark~1.3, one can interpret each subgroup $G$ of $\Cal A$ as a symmetry group
of a spherically homogeneous rooted tree.
Therefore some  concepts related to $G$  were originally stated in terms of  actions on rooted trees.
We restate some of them in an equivalent but more convenient for us standard language of dynamical systems.

For each $g\in G$, we let Supp\,$g:=\{x\in G\mid gx\ne x\}$.
Then Supp\,$g$ is an open subset of $X$.
Given an open subset $O\subset X$, we let
$G_O:=\{g\in G\mid \text{Supp\,}g\subset O\}$.
The group $G_O$ is called {\it the rigid stabilizer of $O$}.
Of course, if $O_1$ and $O_2$ are disjoint open subsets of $X$ then $G_{O_1}$ and $G_{O_2}$ commute.

\definition{Definition 1.7 }\cite{Gr2} Let $G$ be a minimal subgroup of $\Cal A$.
\roster
\item"(i)"  $G$ is called {\it branch}  if  for each $n>0$,
the group
$
\prod_{y\in X_1^n}G_{[y]_1^n}
$
is of finite index in $G$.
\item"(ii)"
$G$ is called {\it weakly branch} if $G_{[y]_1^n}$ is nontrivial for each
 $y\in X_1^n$ and $n\in\Bbb N$.
\endroster
\enddefinition

Experts in the field of rooted trees denote the group $
\prod_{y\in X_1^n}G_{[y]_1^n}
$ by rist$_G(n)$.
Since $G$ is minimal, $G_{[y]_1^n}$ and $G_{[z]_1^n}$ are conjugate in $G$ for all $y,z\in X_1^n$ and $n\in\Bbb N$.
It follows that each branch group is weakly branch.

There is another useful interpretation of the minimal $G$-actions on $X$ (with $G$ being a subgroup in $\Cal A$) as odometer actions (see \cite{DaLe} and \cite{DaVe}).
Fix an element $x=(x_n)_{n=1}^\infty\in X$.
Let
$$
\Gamma_n=\{g\in G\mid \pi_n(gx)=\pi_n(x)\}.
$$
Then $\Gamma_n$ is a cofinite subgroup in $G$ and $\Gamma_1\supset \Gamma_2\supset\cdots$.
We note that  $\Gamma_n$ depends on the choice of $x$.
If we choose another point then the corresponding group will be conjugate to $\Gamma_n$.
Let $\widetilde\Gamma_n=\bigcap_{g\in G}g\Gamma_n g^{-1}$.
In other words, $\widetilde\Gamma_n$ is the largest  subgroup of $\Gamma_n$ that is normal in $\Gamma$.
According to \cite{Gr3, Definition~2.1}, $\widetilde\Gamma_n$ is the $n$-the level stabilizer for the underlying rooted tree.
Of course, $\widetilde\Gamma_n$ is of finite index in $G$ and
$\widetilde\Gamma_1\supset \widetilde \Gamma_2\supset\cdots$.
If $g\in \widetilde\Gamma_n$ and \thetag{1-1} holds for $g$
then $\gamma=I$, $\alpha_1\equiv I$, \dots, $\alpha_{n-1}\equiv I$.
This follows from the minimality of $G$.
This implies that $\bigcap_{n=1}^\infty\widetilde\Gamma_n=\{1\}$.
Hence, {\it the odometer action of $G$ associated with $(\Gamma_n)_{n=1}^\infty$} is well defined
as the projective limit of the following sequence of finite  $G$-spaces:
$$
G/\Gamma_1@<<<G/\Gamma_2@<<<G/\Gamma_3@<<<.
$$
It was shown in \cite{DaVe} that this action is isomorphic to the $G$-action on $X$.

\head 2. Nonsingular group actions
\endhead

In this section we briefly remind some basic concepts of nonsingular dynamics.

Let $(Y,\nu)$ be a non-atomic  probability space.
Denote by Aut$(Y,\nu)$ the group of $\nu$-nonsingular transformations.
The {\it Koopman}  unitary representation $\kappa_\nu$ of Aut$(Y,\nu)$ in $L^2(Y,\nu)$ is defined as follows:
$
\text{Aut}(Y,\nu)\ni g\mapsto \kappa_\nu(g)\in\Cal U(L^2(Y,\nu)),
$
where
$$
\kappa_\nu(g)f:=f\circ g^{-1}\sqrt{\frac{d(\nu\circ g^{-1})}{d\nu}}, \quad f\in L^2(Y,\nu).
$$
Let the unitary group $\Cal U(L^2(Y,\nu))$ is furnished with the weak operator topology.
We endow Aut$(Y,\nu)$ with the weakest topology in which the Koopman representation is continuous.
This topology on Aut$(Y,\nu)$ is called {\it weak}.
Then Aut$(Y,\nu)$ endowed with the weak topology   is a Polish group.

Let $G$ be a countable subgroup of Aut$(Y,\nu)$.
We endow it with the discrete topology.
Then $G$ is called
\roster
\item"(i)" {\it conservative} if for each subset $A\subset Y$ of positive measure and a finite subset $F\subset G$, there is $g\not\in F$ such that
$\nu(gA\cap A)>0$;
\item"(ii)" {\it dissipative} if it is not conservative;
\item"(iii)" {\it rigid} if there is a sequence $(g_n)_{n=1}^\infty$ such that $g_n\to\infty$  and $\kappa_\nu(g_n)\to I$ weakly as $n\to\infty$;
\item"(iv)" {\it ergodic} if each $G$-invariant subset of $Y$ is either $\nu$-null or $\nu$-conull.
\endroster

We note that if $G$ preserves $\nu$ then $G$ is conservative.
If $G$ is rigid then  $G$ is conservative.
If $G$ is ergodic then $G$ is conservative.
There exist non-rigid ergodic transformation groups.
The following criterion of conservativity is well known:
\roster
\item"---"
$G$ is conservative if and only if $\sum_{g\in G}\frac{d(\nu\circ g)}{d\nu}(y)=\infty$ at $\nu$-a.e. $y$.
\endroster

The restriction of $\kappa_\nu$ to $G$ is called the {\it Koopman unitary representation of $G$}.
Properties of $G$ that depends indeed on  the  Koopman unitary representation of $G$
are called {\it spectral}.
Neither conservativity nor ergodicity are spectral properties of $G$.
However, ergodicity is a spectral property within the subclass of $\nu$-preserving transformation groups.
It is easy to see that if $G$ is not ergodic then the Koopman representation $(\kappa_\mu(g))_{g\in G}$ is not irreducible.
However, $(\kappa_\mu(g))_{g\in G}$ can be reducible for ergodic $G$-actions: for instance, if $G$ is Abelian.

Let $X$ be the space of infinite paths on a spherically homogeneous rooted tree and let $\Cal A\subset\text{Homeo\,}(X)$ be as in Section~1.

\proclaim{Proposition 2.1} Let $G\subset \Cal A$ be weakly branch.
Let $\mu$ be a $G$-quasiinvariant Borel  probability on $X$.
%Suppose that $\mu(O)>0$ for each open subset in $X$.
Then $G$ is rigid.
Hence $G$ is conservative.
\endproclaim
\demo{Proof}
We consider first the case where $\mu$ is non-atomic.
Choose a sequence of cylinders $O_1\supset O_2\supset\cdots$
with $\bigcap_n O_n=\{a\}$ for some point $a\in X$.
Hence $\mu(O_n)\to 0$ as $n\to\infty$.
Then we can find a sequence $(g_n)_{n=1}^\infty$ of $G$-elements
such that $\text{Supp\,}g_n\subset O_n$ and  $g_n\ne 1$ for each $n$.
Without loss of generality we may assume (passing to a subsequence, if necessary) that
%$\mu(\text{Supp\,}g_n)> \mu(\text{Supp\,}g_{n+1})$.
 $g_n\ne g_m$ if $n\ne m$.
As $G$ is countable and discrete, we obtain that $g_n\to\infty$.
Take $f,r\in L^2(X,\mu)$.
Then
$$
\align
\langle\kappa_\mu(g_n)f,r\rangle&=
\langle\kappa_\mu(g_n)(f\cdot1_{O_n}),r\rangle+
\langle f\cdot1_{O_n^c},r\rangle\\
&=\langle\kappa_\mu(g_n)(f\cdot1_{O_n}),r\rangle+
\langle f,r\rangle
-\langle f\cdot1_{O_n},r\rangle
\endalign
$$
Since $|\langle\kappa_\mu(g_n)(f\cdot1_{O_n}),r\rangle|\le \|f\|_2\cdot \|1_{O_n}\|_2\cdot \|r\|_2\to 0$
and $|\langle f\cdot1_{O_n},r\rangle|\to 0$ as $n\to\infty$,
it follows that  $\kappa_\mu(g_n)\to I$ weakly as $n\to\infty$.
Thus, $G$ is rigid.

Now consider the case of purely atomic measure.
Let $\mu(x)>0$ for some $x\in X$.
Since $G$ is weakly branch, the stabilizer $H$ of $x$ is an infinite subgroup in $G$.
Select  a sequence $(h_n)_{n=1}^\infty$ of mutually different elements of $H$.
Then $h_n\to\infty$ and $\kappa_\mu(h_n)=I$ for each $n$.
Hence, $G$ is rigid.

The general case reduces  either to the nonatomic one or to the purely atomic one (the weakly branch condition is used essentially for that).
\qed
\enddemo

It follows, in particular,  that the combinatorial property of $G$ to be weakly branch implies, in a sense, a ``universal rigidity'', i.e. $G$ is $\mu$-rigid (and hence conservative) with respect to every(!) quasiinvariant
%nonatomic
measure $\mu$ on $X$.
% such that $\mu(O)>0$ for each open subset $O\subset X$.

\comment

We need a stronger concept than weakly branch.
Let $O$ be a cylinder in $X$ and let
$G_O:=\{g\in G\mid \text{Supp\,}g\subset O\}$.
Then $G$ is weakly branch if and only if $G_O$ is nontrivial for each cylinder $O$.

\definition{Definition 1.6} We call $G$ {\it saturated} if  for each cylinder $O\subset X$,
the group $G_O$ is minimal  on  $O$.
\enddefinition

\endcomment

\head 3. %Dynamical properties of the infinite product measures on  $X$
Finitary groups of automorphisms~of~$T$ for infinite product measures on $\partial T$
\endhead

In this section we find conditions under which groups of automorphisms of $T$ are nonsingular with respect to infinite product measures on $\partial T$.

%\subhead 3.1. Finitary sugroups
%\endsubhead
Let $(X_n)_{n=1}^\infty$, $X$, $\Cal A$ and $\Cal A_0$ be the same  as in Section~1.
For each $n>0$, fix
 a non-degenerared probability $\mu_n$ on $X_n$.
Let
$\mu:=\bigotimes_{n=1}^\infty\mu_n$.
Such probability measures are called infinite product measures on $X$.
Assume that $\mu$ is non-atomic or, equivalently,  $\prod_{n=1}^\infty\max_{x_n\in X_n}\mu_n(x_n)=0$.
From now on we will consider only the nonatomic infinite product measures on $X$.

For $m>n>0$, we let $\mu_n^m:=\mu_n\otimes\cdots\otimes \mu_m$.
%By the Kolmogorov 0-1 law, the tail $\sigma$-algebra on $X$ is trivial (mod $\mu$).
%If each $\mu_n$ is the equi-distribution then we call $\mu$ {\it the Haar measure} on $X$.
%Denote by Aut$(X,\mu)$ the group of $\mu$-nonsingular transformations of $X$.

If $X_1=X_2=\cdots$ and $\mu_1=\mu_2=\cdots$ then the corresponding product measure $\mu$
is often  called {\it Bernoulli}.

\comment

Denote by Aut$(X,\mu)$ the group of $\mu$-nonsingular transformations of $X$.
It is Polish in the weak topology.
We now isolate a subgroup  $\Cal G$ in Aut$(X,\mu)$:
$$
\Cal G:=\{g\in\text{Aut}(X,\mu)\mid g\goth B_n=\goth B_n\text{ for each }n\in\Bbb N\}.
$$
Equivalently, a transformation $g$ from $\text{Aut}(X,\mu)$ belongs to $\Cal G$ if
Then $g$ is continuous (or, more rigorously, coincides with a homeomorphism of $X$ almost everywhere in $\mu$) and preserves the Haar measure on $X$.

\endcomment

\definition{Definition 3.1}
Let $g\in \Cal A$.
Let $\gamma$ and $(\alpha_n)_{n=1}^\infty$ be the objects determined by $g$ as in \thetag{1-1}.

\roster
\item"(i)"
If
 there is $N>0$ such that
$\alpha_n(x_1,\dots,x_n)=I$ for each $n\ge N$ and every $(x_1,\dots,x_n)\in X_1^n$ then $g$ is called  {\it finitary}.
\item"(ii)"
If for every $x=(x_1,x_2,\dots)\in X$, there is $N>0$ such that
$\alpha_n(x_1,\dots,x_n)=I$ for each $n\ge N$  then $g$ is called {\it quasi-finitary}.

\item"(iii)"
If for  $\mu$-a.e. $x\in X$, there is $N>0$ such that
$\alpha_n(x_1,\dots,x_n)=I$ for each $n\ge N$ then $g$ is called {\it $\mu$-finitary}.

\item"(iv)"
If for each $\epsilon>0$, there are $N>0$ and  a clopen subset $A\subset X$ such that
$\mu(A)<\epsilon$ and $g$ changes no more than $N$ first coordinates  of each $x\not\in A$
 then $g$ is called {\it purely $\mu$-finitary}.

\item"(v)"
If for  $\mu$-a.e. $x\in X$, there is $N>0$ such that the bijection
$\alpha_n(x_1,\dots,x_n)$ of $X_{n+1}$ preserves $\mu_{n+1}$ for each $n\ge N$ then $g$ is called {\it eventually $\mu$-presereving}.
\endroster
\enddefinition

It is easy to see that the set of all finitary elements of $\Cal A$ is exactly  $\Cal A_0$.
Denote by $\Cal A_0^+$, $\Cal A_\mu$, $\Cal A_\mu^\bullet$  and $\Cal A_\mu^+$ the set of  quasi-finitary transformations,  the set of $\mu$-finitary transformations, the set of purely $\mu$-finitary transformations
and the set of
eventually $\mu$-preserving transformations respectively.
Then we have:
$$
\Cal A_0\subset\Cal A_\mu^\bullet\subset \Cal A_\mu\quad\text{and}\quad\Cal A_0\subset\Cal A_0^+\subset\Cal A_\mu\subset \Cal A_\mu^+.
$$
A subgroup $G$ of  $\Cal A$ is called {\it finitary, quasi-finitary, $\mu$-finitary, purely $\mu$-finitary or eventually $\mu$-preserving} if
$G\subset \Cal A_0$,  $G\subset \Cal A_0^+$, $G\subset \Cal A_\mu$,
$G\subset \Cal A_\mu^\bullet$
or  $G\subset \Cal A_\mu^+$
respectively.

If a transformation $g$ is $\mu$-finitary then for each $\epsilon>0$ there exist $N>0$ and a Borel subset $A\subset X$ such that $\mu(A)<\epsilon$ and $g$ changes no more than $N$ first  coordinates of each $x\not\in A$.
However, in contrast with the purely $\mu$-finitary transformations, $A$ needs not to be open.
See, for instance, Example~4.2 below.

It is important to note that if a transformation $g$ is eventually $\mu$-presereving then $g$ is $\mu$-nonsingular, i.e. $g\in\text{Aut}(X,\mu)$ and
$$
\aligned
\frac{d\mu\circ g}{d\mu}(x_1,x_2,\dots) &=\lim_{n\to\infty}\frac{d\mu^n_1\circ g_n}{d\mu^n_1}(x_1, \dots,x_n)\\
&=\frac{d\mu^N_1\circ g_N}{d\mu^N_1}(x_1, \dots,x_N)\\
&=\frac{\mu_1^N([g_N(x_1,\dots,x_N)]_1^N)}{\mu_1^N([(x_1,\dots,x_N)]_1^N)}\\
&=
\frac{\mu_1(\gamma x_1)}{\mu_1(x_1)}
\prod_{j=2}^N\frac{\mu_j(\alpha_{j-1}(x_1,\dots,x_{j-1})[x_j])}{\mu_j(x_j)}
\endaligned
\tag3-1
$$
at $\mu$-a.e. $(x_n)_{n=1}^\infty\in X$, where $g_n\in\text{Bij}(X_1^n)$ is determined by $\pi_n\circ g=g_n\circ\pi_n$ and $N$ is defined in Definition~3.1(v).

It is obvious that
$\Cal A_0^+$ is a group.
Since each element of  $\Cal A_\mu$ preserves the ideal of  $\mu$-zero subsets of $X$,
it is easy to verify that $\Cal A_\mu$ is also a group.
In a similar way, one can show that $\Cal A_\mu^+$ is a group too.
We remind that if  a transformation $g$ of $X$ is $\mu$-nonsingular then for each $\epsilon>0$,
there is $\delta>0$ such that if $\mu(A)>\delta$ for a subset $A\subset X$ then $\mu(gA)>1-\epsilon$.
This fact plus a simple observation that each element of $\Cal A$ preserves the ring of clopen subsets in $X$ yield that  $\Cal A_\mu^\bullet$ is a group.
Thus, $\Cal A_0^+$, $\Cal A_\mu$,  $\Cal A_\mu^\bullet$, $\Cal A_\mu^+$
 are  subgroups of $\Cal A$.

The following proposition provides sufficient conditions for a subgroup of $\Cal A$ to be
$\mu$-finitary, eventually $\mu$-preserving or purely $\mu$-finitary.

\proclaim{Proposition 3.2} Let $G$ be a subgroup of $\Cal A$.
Given $n>0$ and $g\in G$, let
$$
\align
F_n(g)&:=\{y\in X_1^n\mid \alpha_n(y)\ne I\}\\
F_n^+(g)&:=\{y\in X_1^n\mid \mu_{n+1}\circ\alpha_n(y)\ne \mu_{n+1}\}\text{ and }\\
F_n^\bullet(g)&:=\{y\in X_1^n\mid\text{there exist $m>n$ and $a\in X_{n+1}^{m}$ with }\alpha_{m}(ya)\ne I\}.
\endalign
$$
\roster
\item"(i)"
If
$
\sum_{n=1}^\infty\mu_1^n(F_n(g))<\infty$
for each $g\in G$
%\tag3-2
then $G$ is $\mu$-finitary.
\item"(ii)"
If
$
\sum_{n=1}^\infty\mu_1^n(F_n^+(g))<\infty$
for each $g\in G$
%\tag3-2
then $G$ is eventually $\mu$-preserving.
\item"(iii)"
$\lim_{m\to\infty}\mu_1^n(F_n^\bullet(g))=0$ if and only if
 $g$ is purely $\mu$-finitary.
\endroster
\endproclaim

\demo{Proof}
(i) and (ii)  follow  from the Borel-Cantelli lemma.
Of course, (iii) holds.
\qed
\enddemo

\remark{Remark 3.3} Let $\Cal R$ stand for the {\it tail} equivalence relation on $X$.
Two points $(x_n)_{n=1}^\infty$ and $(y_n)_{n=1}^\infty$ are called tail equivalent
if $x_n=y_n$ eventually in $n$.
Denote by $[\Cal R]$ the full group of $\Cal R$.
We remind that
$$
[\Cal R]=\{\theta\in\text{Aut}(X,\mu)\mid (x,\theta x)\in\Cal R\text{ at a.e. $x$}\}.
$$
Since $\mu$ is an infinite product measure,
$[\Cal R]$ is ergodic.
Of course,  a subgroup $G$ of $\Cal A$ is $\mu$-finitary if and only if  $(gx,x)\in \Cal R$ at $\mu$-a.e. $x$ for each $g\in G$.
In other words, $G$ is $\mu$-finitary if and only if  $G\subset[\Cal R]$.
Then there  is $n\in\Bbb N\sqcup\{\infty\}$ such that
$\Cal R(x)$ consists of $n$ different $G$-orbits for a.e. $x$ \cite{FeSuZi}.
This $n$ is called the {\it index} of the $G$-orbit equivalence relation in $\Cal R$.
Thus, if $G\subset\Cal A_\mu$  then a natural question arises: what is the index of
the $G$-orbit equivalence relation in $\Cal R$?
\endremark

We now provide several examples illustrating the concepts introduced in Definition~3.1 and Remark~3.3.

\example{Example 3.4} Let $X_1=X_2=\cdots=\{0,1\}$.
Denote by $\theta$ the transposition of $0$ and $1$.
The Grigorchuk group $G$  \cite{Gr1} is generated by four elements: $a,b,c,d$.
These elements are determined by the formulae:
$$
\align
a(x_1,x_2,\dots)&=(\theta x_1,x_2,\dots),\\
b(1^n,0,x_{n+2},\dots)&=
\cases
(1^n,0,\theta x_{n+2},x_{n+3},\dots) &\text{if $n\not\in 3\Bbb Z_+$},\\
(1^n,0,x_{n+2},x_{n+3},\dots)&\text{otherwise}
\endcases
,
\\
c(1^n,0,x_{n+2},\dots)&=
\cases
(1^n,0,\theta x_{n+2},x_{n+3},\dots) &\text{if $n\not\in 1+3\Bbb Z_+$},\\
(1^n,0,x_{n+2},x_{n+3},\dots)&\text{otherwise}
\endcases
,
\\
d(1^n,0,x_{n+2},\dots)&=
\cases
(1^n,0,\theta x_{n+2},x_{n+3},\dots) &\text{if $n\not\in 2+3\Bbb Z_+$},\\
(1^n,0,x_{n+2},x_{n+3},\dots)&\text{otherwise}
\endcases
\endalign
$$
and $b(1,1,\dots)=c(1,1,\dots)=d(1,1,\dots)=(1,1,\dots)$.
%It is easy to verify that $G\subset \Cal A\setminus\Cal A_0^+$, i.e. $G$ is not quasifinitary.
%Consider, for instance, $b$.
%Then
\comment
$$
f_n(b)=
\cases
\{1^{n-1}0\} &\text{if $n-1\not\in 3\Bbb Z_+$},\\
X_1^n &\text{otherwise.}
\endcases
$$
\endcomment
Then $G$ is quasi-finitary but not finitary.
Hence $G$ is $\mu$-nonsingular for each infinite product measure $\mu$ on $X$.
It is straightforward to verify that $G$ is minimal and the $G$-orbit equivalence relation is the tail equivalence relation on $X$.
It follows that $G$ is ergodic with respect to every infinite product measure on $X$.

\endexample

Below is a standard example of a minimal $G$ whose orbit equivalence relation is of index $2$ in the tail equivalence relation.
Examples with any other  index are constructed in a similar way.

\example{Example 3.5} Let $X_1=X_2=\cdots=\Bbb Z/2\Bbb Z$.
Then $X=(\Bbb Z/2\Bbb Z)^{\Bbb N}$.
We consider $X$ as a compact Abelian group.
Denote by $\mu$  the normalized Haar measure on $X$.
Let $\Gamma:=\bigoplus_{n\in\Bbb N}\Bbb Z/2\Bbb Z$.
Then $\Gamma$ is a dense countable subgroup in $X$.
We consider the natural action of  $\Gamma$ on $(X,\mu)$ by translations.
This action is free and preserves $\mu$.
Of course, $\Gamma\subset\Cal A_0$.
We let
$$
G:=\bigg\{(a_n)_{n=1}^\infty\in \Gamma\mid \sum_{n\in\Bbb N}a_n=0\bigg\}.
$$
The sum is well defined since only finitely many coordinates of $(a_n)_{n=1}^\infty$ are non-zero.
Of course, $G$ is a subgroup of index 2 in $\Gamma$.
It is straightforward to verify that $G$ is minimal.
Let $\Cal R$ and $\Cal S$ stand for $\Gamma$- and $G$-orbit equivalence relations respectively.
Then $\Cal R$ is the tail equivalence relation on $X$ and $\Cal S$ is a subrelation of $\Cal R$ of index 2.
Thus, $G$ is minimal, finitary and the $G$-orbit equivalence relation is of index 2 in the tail equivalence relation on $X$.
\endexample

Given two orthogonal measures on $X$, we can construct a minimal $G$ which is
 finitary with respect to one of these measures but not finitary with respect to the other one.

\example{Example 3.6}
Let $X$ be as in Example~3.5 and
let $\mu$ and $\nu$ be two arbitrary nonatomic  probabilities on $X$ such that $\mu\perp\nu$ and
$\mu(O)>0$ and $\nu(O)>0$ for each open subset $O\subset X$.
We do not even assume that they are of infinite product type.
Then there exists a Borel subset $B\subset X$ such that $\mu(B)=1$ and $\nu(B)=0$.
Given $\epsilon>0$, we can find an open subset $O\subset X$ such that $O\supset B$ and $\nu(O)<\epsilon$.
Of course, $\mu(O)=1$.
Since $O=\bigcup_{n=1}^\infty O_n$  for a sequence $O_1\subset O_2\subset\cdots$ of clopen subsets $O_n\subset X$,
we can find $n$ such that $\mu(O_n)>1-\epsilon$.
Of course, $\nu(O_n)<\epsilon$.
Since $O_n$ is clopen, there is $m>0$ and a subset $D\subset X_1^m$ such that
$O_n=\bigsqcup_{d\in D}[d]_1^m$.

Using this reasoning repeatedly, we can construct an increasing sequence $n_1<n_2<\cdots$ of positive integers  and a sequence of  subsets $D_k\subset X_1^{n_k}$ such that if
$U_k=\bigsqcup_{d\in D_k}[d]_1^{n_k}$ then $\mu(U_k)>1-2^{-k}$ and $\nu(U_k)<2^{-k}$ for each $k\in\Bbb N$.
Let $\tau_k$ be the non-trivial bijection of $X_{n_k+1}$, $k=1,2,\dots$.
Let $\gamma$ be the trivial bijection of $X_1$.
We now define  a sequence of mappings
$\alpha_m:X_{1}^m\to\text{Bij}(X_{m+1})$ by setting
$$
\alpha_{n_k}(y)=
\cases
\tau_{k} &\text{if $y\in D_k$},\\
I &\text{if $y\not\in D_k$}
\endcases
$$
and
$\alpha_m\equiv I$ if $m\ne n_k$ for any $k$.
Let $g\in\Cal A$ be a homeomorphism of $X$ that is determined by $\gamma$ and $(\alpha_n)_{n=1}^\infty$ via \thetag{1-1}.
Let $F_n(g)$ stand for the subset defined in the statement  Proposition~3.2, $n\in\Bbb N$.
Then $F_m(g)=D_k$ if $m=n_k$ and $F_m(g)=\emptyset$ otherwise.
It follows from Proposition~3.2(i) that
 $g\in \Cal A_\nu$.
The Borel-Cantelli lemma implies also that
 $g\notin \Cal A_\mu$.
Let $G$ be generated by $\Cal A_0$ and $g$.
Then $G$ is a minimal transformation group, $G$ is $\nu$-finitary but not $\mu$-finitary.
\endexample

In the following example, for each Bernoulli measure on $\{0,1\}^\Bbb N$, we construct a minimal $G$-action that is finitary with respect to this measure and non-finitary with respect to any other Bernoulli measure on this space.

\example{Example 3.7} Let $X_1=X_2=\cdots=\{0,1\}$.
For each $\lambda\in(0,1)$, let $\kappa_\lambda$ be the probability on $\{0,1\}$ such that
$\kappa_\lambda(1)=\lambda$.
Let $\mu_\lambda:=(\kappa_\lambda)^{\otimes\Bbb N}$ denote the corresponding
Bernoulli measure.
Fix $\theta\in(0,1)$.
We are going to construct a minimal $G$ which is $\mu_\theta$-finitary but not $\mu_\lambda$-finitary for each $\lambda\ne\theta$.
For that we consider the probability space $(X,\mu_\lambda)$ and  random variables $f_n:X\to\Bbb R$ given by $f_n(x)=x_n$, where $x_n$ is the $n$-th coordinate of $x$,  $n\in\Bbb N$.
Then $(f_n)_{n\in\Bbb N}$ are independent and equally distributed.
By the weak law of large numbers,
$$
%\mu_\lambda(\{x\in X\mid n^{-1}(x_1+\cdots+x_n)\to \lambda\})\to 1
n^{-1}(f_1+\cdots+f_n)\to \lambda\quad\text{in $\mu_\lambda$}
$$
as $n\to\infty$.
We choose a decreasing sequence $\epsilon_n\to 0$ of reals and a sequence of closed subsets
$\Lambda_1\subset\Lambda_2\subset\cdots$ of $(0,1)$ such that
$$
\gather
\bigcup_{n=1}^\infty\Lambda_n=(0,\theta)\cup(\theta,1)\quad
\text{ and }\tag3-2\\
\min_{\lambda\in\Lambda_n}|\theta-\lambda|>5\epsilon_n\quad\text{ for each $n$.}\tag3-3
\endgather
$$
For $\lambda\in(0,1)$, $n\in\Bbb N$ and $\epsilon>0$, we let
$$
B_{\lambda,n,\epsilon}:=\{x\in X_1^n\mid |n^{-1}(x_1+\cdots+x_n)-\lambda|<\epsilon\}.
$$
It follows from the weak  law of large numbers that for each $\eta>0$,
$$
(\kappa_\lambda)^{\otimes n}(B_{\lambda,n,\epsilon})>1-\eta\tag3-4
$$
eventually in $n$.
We note that if \thetag{3-4} holds for some $\lambda$ then \thetag{3-4} holds also for a
   neighborhood of $\lambda$ (provided that $n$, $\epsilon$ and $\eta$ are fixed).
Therefore, using \thetag{3-4} and the compactness of $\Lambda_k$, we can find a sequence $n_1<n_2<\cdots$ of positive reals such that for each $k$,
$$
(\kappa_\theta)^{\otimes n_k}(B_{\theta,n_k,\epsilon_k})>1-k^{-2}\quad\text{and}\quad
\min_{\lambda\in\Lambda_k}(\kappa_\lambda)^{\otimes n_k}(B_{\lambda,n_k,\epsilon_k})>1-k^{-2}.
$$
It follows from \thetag{3-3} that $B_{\theta,n_k,\epsilon_k}\cap B_{\lambda,n_k,\epsilon_k}=\emptyset$ for each $\lambda\in \Lambda_k$.
Let $\gamma$ be the trivial bijection of $X_1$.
Let $\tau_k$ be the non-trivial bijection of $X_{n_k+1}$.
We now define a sequence of mappings
$\alpha_m:X_{1}^m\to\text{Bij}(X_{m+1})$ by setting
$$
\alpha_{n_k}(y)=
\cases
\tau_{k} &\text{if $y\not\in B_{\theta,n_k,\epsilon_k}$},\\
I &\text{if $y\in B_{\theta,n_k,\epsilon_k}$}
\endcases
$$
and
$\alpha_m\equiv I$ if $m\ne n_k$ for any $k$.
Let $g\in\Cal A$ be a homeomorphism of $X$ that is determined by $\gamma$ and $(\alpha_n)_{n=1}^\infty$ via \thetag{1-1}.
Let $F_n(g)$ stand for the subset defined in the statement  Proposition~3.2, $n\in\Bbb N$.
Then
$$
F_m(g)=
\cases
(B_{\theta, n_k,\epsilon_k})^c&\text{if $m=n_k$ for some $k$,}\\
\emptyset &\text{otherwise.}
\endcases
$$
It follows that $\sum_{m=1}^\infty (\kappa_\theta)^{\otimes m}(F_{m}(g))<\sum_{k=1}^\infty k^{-2}<\infty$.
Hence, by Proposition~3.2(i),
$g\in \Cal A_{\mu_\theta}$.
On the other hand,
$$
\align
(\kappa_\lambda)^{\otimes n_k}((B_{\theta,n_k,\epsilon_k})^c)
&=
1-(\kappa_\lambda)^{\otimes n_k}(B_{\theta,n_k,\epsilon_k})\\
&\ge 1-
(\kappa_\lambda)^{\otimes n_k}((B_{\lambda,n_k,\epsilon_k})^c)\\
&>
(\kappa_\lambda)^{\otimes n_k}(B_{\lambda,n_k,\epsilon_k})\\
&>1-k^{-2}.
\endalign
$$
Therefore,
$\sum_{m=1}^\infty (\kappa_\lambda)^{\otimes m}(F_{m}(g)^c)<\infty$.
It follows from the Borel-Cantelli lemma that
 $g\notin \Cal A_{\mu_\lambda}$
for each $\lambda\in\Lambda_{k}$.
In view of \thetag{3-2}, $g\notin \Cal A_{\mu_\lambda}$ for each $\lambda\ne\theta$.
Let $G$ is generated by $\Cal A_0$ and $g$.
Then $G$ is a minimal transformation group, $G$ is $\mu_\theta$-finitary
but
  not $\mu_\lambda$-finitary if $\lambda\ne\theta$.

\endexample

%where $g_N$ is the bijection of $X^N_1$ given by
%by the formula
 %$$
% g_N\pi_N(x)=\pi_n(gx)\quad \text{ for a.e. $x\in X$.}
% $$
%Denote by $G_n$ the subgroup in Bij$(X_1^n)$ generated by  $g_n$ when $g$ runs $G$.
%It follows from the remark:
%\roster
%\item"---"
% $G$ is ergodic with respect to the Haar measure if and only if $G_n$ is transitive for each $n$.
% Moreover, $G$ is uniquely ergodic.

% \item"---"
% If $G$ is $\mu$-finitary and ergodic with respect to $\mu$ then $G_n$ is transitive
%for each $n$.
% \endroster

\head 4. Subexponential boundedness and ergodicity
\endhead

\subhead 4.1.  On subexponential boundedness
\endsubhead
 In this subsection we compare the condition from Proposition~3.2(i) with the concept of
subexponential boundedness that was introduced in  \cite{DuGr1} in the case where
$X_1=X_2=\cdots$.
Let $g\in \Cal A$.
We remind that $F_n(g)$ and $F_n(g)^\bullet$ were defined in Proposition~3.2.

\definition{Definition 4.1 } Let $X_1=X_2=X_3=\cdots$.
\roster
\item"(i)"
If   $\sum_{n=1}^\infty \#(F_n^\bullet(g))\delta^n<\infty$ for each $\delta\in(0,1)$
and  $g\in G$ then $G$ is called {\it subexponentially bounded} \cite{DuGr1}\footnote{A.Dudko and R.Grigorchuk   used the  condition $ \lim_{n\to\infty}\#(F_n^\bullet(g))\delta^n= 0$   for each $\delta\in(0,1)$ in the definition of subexponential boundedness in \cite{DuGr1}. Of course, this condition  is equivalent to Denition~4.1(i).}.
\item"(ii)"
If   $\sum_{n=1}^\infty \#(F_n(g))\delta^n<\infty$ for each $\delta\in(0,1)$
and  $g\in G$ then we call $G$  {\it (w)-subexponentially bounded}.
\endroster
\enddefinition

Since $F_n^\bullet(g)\subset  F_n(g)$,
every  subexponentially bounded group is
(w)-subexponen\-tially bounded.
We note also that $\# F_n(g)\le(\# X_1)^n$ for each $g\in G$ and $n\in\Bbb N$.
Hence
$$
\sum_{n=1}^\infty\# (F_n(g))\delta^n<\infty\quad\text{whenever $0<\delta<\frac 1{\# X_1}$.}
$$
It was shown in \cite{DuGr1} that if $G$ is subexponentially bounded  and $\mu$ is  Bernoulli
then $G$ is purely $\mu$-finitary.
In particular, $G$ is $\mu$-nonsingular.
We  note also
that each  (w)-subexponentially bounded group
is $\mu$-finitary.
This
follows from Proposition~3.2(i) because
$$
\mu_1^n(F_n(g))\le \#(F_n(g))\max_{y\in X_1^n} \mu_1^n(y)=\#(F_n(g))\Big(\max_{x_1\in X_1}\mu_1(x_1)\Big)^n
$$
and, hence, the inequality from Proposition~3.2(i) holds.

Thus, the (w)-subexponential boundedness implies a sort of ``universal nonsingularity'', i.e. if $G$ is
(w)-subexponentially bounded then $G$
is nonsingular for every Bernoulli measure on $X$.

In the following two examples we show that the concept of
(w)-subexponential boundedness is more general than  the subexponential boundness and that the
$\mu$-finitarity is  more general
than the (w)-subexponential boundedness.

\example{Example 4.2} Let $X_1=X_2=\cdots=\{0,1\}$.
For each $n>0$, we define a mapping $\alpha_n:X_1^n\to\text{Bij}(X_1)$ in the following way.
We set $\alpha_n(y)=I$  for every $y\in X_1^n$ if $n\ne k!$ for any $k>1$ and
$$
\alpha_{k!}(y_1,\dots,y_{k!})\ne I\iff\text{ $y_{(k-1)!+1}=y_{(k-1)!+2}=\cdots=y_{k!}=0$}.
$$
Then we define  a homeomorphism $g\in\Cal A$ by setting
$$
g(x_1,x_2,\dots)=(x_1,\alpha_1(x_1)[x_2], \alpha_2(x_1,x_2)[x_3],\dots).
$$
Let $G$ be the cyclic group generated by $g$.

If $n\in\Bbb N$ and $y\in X_1^n$ then we choose $k$ such that $(k-1)!>n$.
Take $x=(x_j)_{j=1}^{k!}\in X_1^{k!}$ such that $x_1=y_1$, \dots, $x_n=y_n$ and $x_{(k-1)!+1}=\cdots=
x_{k!}=0$.
Then $\alpha_{k!}(x)\ne I$.
Hence, $y\in F_n^\bullet$.
Thus, $F_n^\bullet=X_1^n$.
Therefore,  if $\delta>0.5$ then
$$
\sum_{n=1}^\infty\# (F_n^\bullet)\delta^n=\sum_{n=1}^\infty(2\delta)^n=\infty.
$$
This means that  $G$ is not subexponentially bounded.
On the other hand,
it is straightforward to verify that
$$
F_n(g)=
\cases
\emptyset  &\text{if $n\not\in\{2!, 3!,\dots\}$},\\
 X_1^{(k-1)!}0^{k!-(k-1)!} &\text{if $n=k!$ for some $k>1$}.
\endcases
$$
Hence,
$$
\sum_{n=2}^\infty\# (F_n)\delta^n=\sum_{k=2}^\infty \# (F_{k!})\delta^{k!}=
\sum_{k=2}^\infty 2^{(k-1)!}\delta^{k!}
=\sum_{k=2}^\infty (2\delta^k)^{(k-1)!}<\infty
$$
for each $\delta\in(0,1)$.
This means that $G$ is (w)-subexponentially bounded.
Hence, $G$ is $\mu$-finitary for every Bernoulli measure $\mu$ on $X$.
We also note that if $\mu$ is an arbitrary  infinite product probability
measure on $X$ then $\mu_1^n(F_n^\bullet)=1\not\to 0$
as $n\to\infty$.
Hence, $g$ is not purely $\mu$-finitary.
\endexample

Prior to come to the next example,
we remind the Shannon-MacMillan theorem.
The {\it entropy} of $\mu_1$ is a real number
$$
h(\mu_1):= -\sum_{x_1\in X_1}\mu_1(x_1)\log\mu_1(x_1).
$$
We have that $0<h(\mu_1)\le\log\#X_1$.
Moreover, $h(\mu_1)=\log\#X_1$ if and only if $\mu_1$ is the equidistribution.

\proclaim{Shannon-McMillan theorem}
Let  $X_1=X_2=\cdots$ and $\mu_1=\mu_2=\cdots$.
For each $\delta>0$, let
$Y_{n,\delta}:=\Big\{y\in X_1^n \,\Big|\,
\,\Big| \frac1n \log\mu_1^n(y)+h(\mu)\Big|\, <\delta\Big\}$.
Then
$
\lim_{n\to\infty}\mu_1^n(Y_{n,\delta})=
1.
%\quad\text{and}\quad\lim_{n\to\infty}\frac{\log \#Y_{n,\delta}}n=h(\mu_1).
$

\endproclaim

\example{Example 4.3}
Let $X_1=X_2=\cdots$ and $\mu_1=\mu_2=\cdots$.
It follows from the Shannon-McMillan theorem that
 there exist  subsets $B_n\subset X_1^n$ such that
 \roster
\item"(i)" $\mu^n_1(B_n)\ge n^{-2}$ for each $n$ and
\item"(ii)" $\mu_1^n(B_n)n^2\to 1$ as $n\to\infty$ and
 \item"(iii)"
 $\sup_{y\in B_n}\Big|\frac1n \log\mu_1^n(y)+h(\mu)\Big|\to0$ as $n\to\infty$.
\endroster
Let $\gamma$ be a transitive cyclic permutation on $X_1$.
Define a mapping $\beta_n:X_1^n\to\text{Bij}(X_1)$ by setting
$$
\beta_n(y)=
\cases
\gamma &\text{if $y\in B_n$,}\\
I &\text{otherwise.}
\endcases
$$
We now define a transformation $g_n$ of $X$ as follows:
$$
g_nx:=(x_1,\dots,x_{n-1},\gamma x_n,\beta_n(x_1,\dots,x_n)[x_{n+1}],
\beta_{n+1}(x_1,\dots,x_{n+1})[x_{n+2}], \dots).
$$
for each $x=(x_1,x_2,\dots)\in X$.
Of course, if $m>n$ then $F_m(g_n)=B_m$.
Hence, $\sum_{m=1}^\infty \mu_1^m(F_m(g_n))<\infty$
in view of (ii).
By Proposition~3.2(i), $g_n\in \Cal A_\mu$.
Let $G$ be generated by all $g_n$, $n\in\Bbb N$.
Then $G$ is $\mu$-finitary.
Of course, $G$ is topologically transitive and hence minimal (see Claim~1.4).
We claim that $e^{-h(\mu)}$ is a ``bifurcation point''
 in the following sense:
 \roster
\item"(a)"  $\sum_{m=1}^\infty \#F_m(g_n)\delta^m<\infty$ if  $0<\delta< e^{-h(\mu)}$ and
\item"(b)"
$\sum_{m=1}^\infty \#F_m(g_n)\delta^m=\infty$ if  $\delta>e^{-h(\mu)}$
\endroster
for each $n$.
Fix $\epsilon>0$ and take $m>n$.
As %$B_n\subset Y_n$,
(iii) holds,
we obtain that
$$
\frac{\mu_1^m(F_m(g_n))}{e^{-m(h-\epsilon)}}
<\#F_m(g_n)
<\frac{\mu_1^m(F_m(g_n))}{e^{-m(h+\epsilon)}}
$$
eventually in $m$.
Hence,
$$
\sum_{m>n}\#F_m(g_n)e^{-m(h-2\epsilon)}\ge
\sum_{m>n}\frac1{m^2}e^{m\epsilon}=+\infty.
$$
This implies (b).
In a similar way one can prove (a).
Hence, $G$ is not (w)-subexponentially bounded.
Thus, $G$ is minimal, $\mu$-finitary  but not (w)-subexponen\-tially bounded.

\endexample

\subhead 4.2. On ergodicity of subexponentially bounded groups with respect to Bernoulli measures
\endsubhead
It is claimed in   \cite{DuGr1, Proposition~3} that if $G$  subexponentially bounded and minimal as a transformation group on $X$ then it is ergodic with respect to each Bernoulli measure on $X$.
In this subsection we provide 3 counterexamples to this claim.

In the first one  we construct for each Bernoulli measure $\mu$ on $\{0,1\}^\Bbb N$ with $\mu_1(0)\ne\mu_1(1)$,
 a subexponentially bounded group $G$ which is minimal but dissipative with respect to $\mu$.

\example{Example 4.4} Let $X_1=X_2=\cdots=\{0,1\}$ and $\mu_1=\mu_2=\cdots$.
Assume that $\mu_1(0)>\mu_1(1)$.
Fix a sequence of positive reals $(\epsilon_n)_{n=1}^\infty$ such that $\sum_{n=1}^\infty\epsilon_n<\infty$.
Denote by $\gamma$ the non-identity bijection of $X_1$.
Given $n>0$ and $x=(x_j)_{j=1}^\infty\in X$, we let
$A_n(x):=\frac1n(x_1+\cdots+x_n)$.
By the weak law of large numbers, the sequence of random variables
$$
A_n:X\ni x\mapsto A_n(x)\in\Bbb R
$$
converges in measure $\mu$ to the constant $\mu_1(1)$  as $n\to\infty$.
Applying this law infinitely many times and using the fact that $\mu_1(1)<0.5$, we can construct (inductively)
an increasing  sequence $1=n_0<n_1<n_2<\cdots$ of integers such that
$$
\mu_{n_k+1}^{n_{k+1}}\Big(\Big\{(x_j)_{j=n_k+1}^{n_{k+1}}\in X_{n_k+1}^{n_{k+1}}\mid
 \frac1{n_{k+1}-n_k}\sum_jx_j>\frac 12\Big\}<\epsilon_k
$$
for each $k\ge 0$.
Let
$$
\align
Y_k&:=\Big\{(x_j)_{j=n_k+1}^{n_{k+1}}\in X_{n_k+1}^{n_{k+1}}\mid
 \frac1{n_{k+1}-n_k}\sum_jx_j<\frac 12\Big\}\quad\text{and}\\
 Y&:=Y_1\times Y_2\times Y_3\times\cdots\subset X.
 \endalign
 $$
 Then
$\mu_{n_k+1}^{n_{k+1}}(Y_k)>1-\epsilon_k$ and hence
$\mu(Y)=\prod_{k=0}^\infty \mu_{n_k+1}^{n_{k+1}}(Y_k)>0$.
For each $x=(x_j)_{j=n_k+1}^{n_{k+1}}\in Y_k$,
we have that
 $$
 \underbrace{\gamma\times\cdots\times\gamma}_{ n_{k+1}-n_k \text{ times}} x=(1-x_j)_{j=n_k+1}^{n_{k+1}}\notin Y_k.
 $$
Hence $ \underbrace{\gamma\times\cdots\times\gamma}_{  n_{k+1}-n_k \text{ times}}Y_k\cap Y_k=\emptyset$ for each $k\ge0$.
We now define a homeomorphism $\gamma_k$ of $X$ by setting $\gamma_kx=((\gamma_kx)_s)_{s=1}^\infty$ for each $x\in X$, where
$$
(\gamma_kx)_s:=
\cases
\gamma x_s &\text{ if either $s=k$ or $n_k<s\le n_{k+1}$},\\
x_s &\text{otherwise.}
\endcases
$$
Of course, $\gamma_k\in\Cal A_0$, $\gamma_k^2=\text{id}$ and $\gamma_k\gamma_l=\gamma_l\gamma_k$ for all $k,l\ge 1$.
Since $\gamma_k$ is finitary, $\gamma_k$ is $\mu$-nonsingular
and subexponentialy bounded.
Denote by $G$ the group generated by $\gamma_k$, $k\in\Bbb N$.
Then $G$ is Abelian.
It is straightforward to verify that $G$ is minimal as a transformation group on $X$.
We claim that $G$ is not conservative on $(X,\mu)$.
For that, we will show that $Y$ is a wandering subset, i.e. $gY\cap Y=\emptyset$ for each $g\in G\setminus\{1_G\}$.
Indeed, if $g\in G$ then there is $k>0$ such that $g=\gamma_1^{i_1}\gamma_2^{i_2}\cdots\gamma_{k-1}^{i_{k-1}}\gamma_k$ for some $i_1,\dots, i_{k-1}\in \{0,1\}$.
Then
$$
gY\subset g(X_1^{n_k}\times Y_k\times Y_{k+1}\times\cdots)=X_1^{n_k}\times
(\underbrace{\gamma\times\cdots\times\gamma}_{  n_{k+1}-n_k \text{ times}})Y_k\times Y_{k+1}\times\cdots.
$$
It follows that $gY\cap Y=\emptyset$, as claimed.
\endexample

In the following example we show how to modify  Example 4.4 to obtain a conservative non-ergodic minimal action of a subexponentially bounded group.

\example{Example 4.5} Let $(X,\mu)$, $\gamma$,  $(n_k)_{k=0}^\infty$ and $(Y_k)_{k=1}^\infty$ be the same as in  Example~4.4.
Let $\widetilde Y_k$ be the projection of $Y_k$, which is a subset of
$X_{n_k+1}^{n_{k+1}}=X\times X_{n_k+2}^{n_{k+1}}$, onto
$X_{n_k+2}^{n_{k+1}}$.
We now let
$$
\widetilde Y:=X_{1}\times \widetilde Y_1\times X_{n_1+1}\times \widetilde Y_2
\times X_{n_2+1}\times \widetilde Y_3\times\cdots.
$$
Then $\mu(\widetilde Y)>0$.
We now define a homeomorphism $\widetilde\gamma_k$ of $X$ by setting $\widetilde\gamma_kx=((\widetilde\gamma_kx)_s)_{s=1}^\infty$ for each $x\in X$, where
$$
(\widetilde\gamma_kx)_s:=
\cases
\gamma x_s &\text{ if either $s=k$ or $n_k+1<s\le n_{k+1}$},\\
x_s &\text{otherwise.}
\endcases
$$
Of course, $\widetilde\gamma_k\in\Cal A_0$, $\widetilde\gamma_k^2=\text{id}$ and $\widetilde\gamma_k\widetilde\gamma_l=\widetilde\gamma_l\widetilde\gamma_k$ for all $k,l\ge 1$.
Since $\widetilde\gamma_k$ is finitary, $\widetilde\gamma_k$ is $\mu$-nonsingular.
Denote by $G$ the group generated by $\widetilde\gamma_k$,
$k\in\Bbb N$.
Then $G$ is Abelian, minimal and
$$
g\widetilde Y\cap \widetilde Y=\emptyset
\quad\text{ for each $g\in G\setminus\{1_G\}$.}
\tag4-1
$$
For each $k>0$, consider a transformation $\delta_k$ of $X$ which changes only the $n_k$-th coordinate.
Of course, it moves the $n_k$-th coordinate exactly the same way as $\gamma$ does.
Then $\delta_k\in\Cal A_0$, $\delta_k^2=\text{id}$ and $\delta_k\delta_l=\delta_l\delta_k$ for all $k,l\ge 1$.
Denote by $H$ the group generated by all $\delta_k$, $k\in\Bbb N$.
Then $H$ is an Abelian group commuting with $G$.
Let $M$ stand for the group generated by $G$ and $H$.
Then $M$ is Abelian and minimal and
subexponentially bounded because it consists of finitary transformations.
We claim that $M$ is $\mu$-conservative.
Indeed, since $\frac{d\mu\circ \delta_k}{d\mu}(x)\ge\frac{\mu_1(1)}{\mu_1(0)}$ at a.e. $x$,
it follows that
$$
\sum_{h\in H}\frac{d\mu\circ h}{d\mu}(x)>
\sum_{k\in \Bbb N}\frac{d\mu\circ \delta_k}{d\mu}(x)=+\infty\quad\text{at a.e. $x\in X$}.
$$
Hence, $H$ is $\mu$-conservative.
As $H$ is a subgroup of $M$, we obtain that $M$ is conservative, as desired.
We note that $\widetilde Y$ is invariant under $H$.
It is straightforward to verify that the restriction of $H$ to $\widetilde Y$ is not ergodic.
Hence there are two disjoint $H$-invariant subsets $A$ and $B$ in $\widetilde Y$ with
$\mu(A)\mu(B)>0$.
%We claim that $M$ is not ergodic.
%Suppose that the contrary holds: $M$ is ergodic.
%Since $\widetilde Y$ is $H$-invariant,
It follows that the two subsets $\bigcup_{g\in G}g A$ and
$\bigcup_{g\in G}g B$ of $X$ are
$M$-invariant and of positive measure.
We deduce from \thetag{4-1} that these subsets are disjoint.
Hence $M$ is not ergodic.
\endexample

In the following example we  `refine' further Example 4.5.
We show that  $M$ can be chosen weakly branch.
Then, by Proposition~2.1, $M$ is rigid,   not just  conservative.

\example{Example 4.6}
Let $(X,\mu)$,  $(n_k)_{k=0}^\infty$, $\widetilde Y$,  $G$, $(\widetilde \gamma_k)_{k=1}^\infty$
and $(\delta_k)_{n=1}^\infty$ be the same as in ~Example~4.5.
For each $k\in\Bbb N$, $y\in X_1^{n_k-1}$ and $x\in X$, we let
$$
\delta_{k,y}x:=
\cases
\delta_kx & \text{if $x\in [y]_1^{n_k-1}$},\\
x&\text{otherwise.}
\endcases
$$
Then $\delta_{k,y}\in\Cal A$, $\delta_{k,y}\delta_{k,y'}=\delta_{k,y'}\delta_{k,y}$ for all $y,y'\in X_1^{n_k-1}$ and $\prod_{y\in X_1^{n_k-1}}\delta_{k,y}=\delta_k$ for each $k$.
Let $H$ be the group generated by $\delta_{k,y}$, $y,\in X_1^{n_k-1}$, $k\in \Bbb N$.
Then $\widetilde Y$ is invariant under $H$.
Let $M$ be the group generated by $G$ and $H$.
Of course, $M\subset\Cal A_0$.
Moreover, $M$ is minimal and weakly branch.
It is straightforward to verify that $H$ is normal in $M$.
Therefore, each element $m\in M$ can be written  as  $m=hg$ for some  $h\in H$ and $g\in G$.
Following the argument in Example~4.5 almost literally, we obtain that
$M$ is not ergodic.
\endexample

This example demonstrates that even the two properties ``weakly branch''  and ``subexponential boundedness''  together (used in \cite{DoGr1, Theorem~3, 1)}) do not imply ergodicity of the underlying group action with respect to Bernoulli measures.
Hence, the corresponding Koopman representation is not irreducible.
%Thus,  the statement of \cite{DoGr1, Theorem~3, 1)} is not true.

\comment

\example{Example 2.5} Let $X_1=X_2=\cdots$ and $\mu_1=\mu_2=\cdots$.
We  denote the cardinality of $X_1$ by $p$.
Assume that $\mu_1$ is not the equidistribution.
Suppose that we have constructed
\roster
\item"---" an increasing  sequence $1=n_0<n_1<n_2<\cdots$ of integers,
\item"---"  subsets
$Y_k\subset X_{n_{k-1}+1}^{n_{k}}$, $k=1,2,\dots$, and
\item"---"  permutations $\gamma_k,\theta_{k}\in\text{Bij}(X^{n_k}_{n_{k-1}+1})$, $k=1,2,\dots$,
\endroster
 such that
 \roster
  \item"(i)"
  $\mu(X_1\times Y_1\times Y_2\times\cdots)>0$ and, for each $k>0$,
 \item"(ii)"
$\gamma_k$ is of order $p^{n_{k-1}-n_{k-2}}$,
\item"(iii)" the subsets $\gamma_k^jY_k$, $0\le j<p^{n_{k-1}-n_{k-2}}$,
are pairwise disjoint,
 \item"(iv)"
$\theta_k$ is of order
  $p^{n_{k}-2n_{k-1}+n_{k-2}}$,
  \item"(v)"
  $\theta_k\gamma_k=\gamma_k\theta_k$,
   \item"(vi)"
   the group generated by $\theta_k$ and $\gamma_{k}$ is transitive on
   $X^{n_k}_{n_{k-1}+1}$,
     \item"(vii)"
the  transformation $g_k$ of $X$, given by
$$
g_kx:=(x_1,\dots, x_{n_{k-2}}, \theta_{k-1}(x_{n_{k-2}+1},\dots,x_{n_{k-1}}),
\gamma_k(x_{n_{k-1}+1},\dots,x_{n_k}), x_{n_{k}+1},\dots)
$$
at each $x=(x_1,x_2,\dots,)\in X$,
belongs to $\Cal A$.
\endroster
%Then $g_k$ has the same order as $\gamma_k$.
We also define a transformation $g_0\in \Cal A$ by setting
$$
g_0(x_1,x_2,x_3,\dots):=(\theta_0x_1,x_2, x_3,\dots),
$$
where $\theta_0$ is a transitive permutation of $X_1$.
Let $G_n$ be the group generated by $g_k$ for all $k\le n$.
Let $G=\bigcup_{n=0}^\infty G_n$.
Then $G$ is an Abelian group.
Moreover,  $G\subset\Cal A$.
Of course, $G$ is minimal, weakly branch and
subexponentially bounded.
We now show that $G$ is not conservative.
Indeed, let $B:=X_1\times Y_1\times Y_2\times\cdots$.
Then $\mu(B)>0$.
We claim that $B$ is a $G$-wandering subset of $X$.
Take $g\in G\setminus\{1_G\}$.
As $G_1\subset G_2\subset\cdots$ with $\bigcup_{k=1}^\infty G_k=G$, there is $k>0$ such that
$g\in G_k\setminus G_{k-1}$.
Then there exist   $j\in \{1,\dots, p^{n_{k-1}-n_{k-2}}-1\}$ and $h\in G_{k-1}$ such that
$g=hg_{k}^j$.
Hence
$$
gB=\cdots \times \gamma_k^j Y_k\times Y_{k+1}\times Y_{k+2}\times\cdots
$$
Since $\gamma_k^j Y_k\cap Y_k=\emptyset$,
it follows that  $gB\cap B=\emptyset$, as desired.

Thus, it remains to construct the sequence $(n_k,Y_k,\theta_k, \gamma_k)_{k=1}^\infty$ satisfying (i)--(vii).
We will construct it inductively.
Fix
 a sequence $(\epsilon_n)_{n=1}^\infty$ of positive reals  such that
$\sum_{n=1}^\infty\epsilon_n<0.5$.

On the first step of the construction, fix a  transitive permutation  $\theta_0$ of $X_1$.
Since $\mu_1$ is not the equidistribution, $h(\mu_1)<\log p$.
Applying the Shannon-McMillan theorem, we find $n_1>0$ and a subset $Y_1\subset X_2^{n_1}$
such that
$$
\mu_2^{n_1}(Y_1)>1-\epsilon_1\quad\text{and}\quad \# Y_1\cdot p<p^{n_1-1}.\tag2-2
$$
Given two points $u:=(x_2,\dots,x_{n_1})$ and $u':=(x_2,\dots,x_{n_1})$ in $X_2^{n_1}$, we define
a permutation  $\psi_{u,u'}$ of $X_2^{n_1}$ by setting
$$
\psi_{u,u'}(y_2,\dots,y_{n_1}):=(\tau y_2, \beta_1(y_2)[y_3],\dots, \beta_{n-1}(y_2,...,y_{n_1-1})[y_{n_1}]),\tag2-3
$$
where $\tau$ is the transposition of $x_2$ and $x_2'$  and
$$
\beta_j(y_2,\dots,y_{k-1})[y_{k}]:=
\cases
y_k, &\text{if $(y_2,\dots,y_{k-1})\ne (x_2,\dots,x_{k-1})$}\\
y_k, &\text{if $(y_2,\dots,y_{k-1})= (x_2,\dots,x_{k-1})$ but $y_k\not\in\{x_k,x_k'\}$}\\
x_k', &\text{if $(y_2,\dots,y_{k-1})= (x_2,\dots,x_{k-1})$ and  $y_k=x_k$}\\
x_k, &\text{if $(y_2,\dots,y_{k-1})= (x_2,\dots,x_{k-1})$  and  $y_k=x_k'$}
\endcases
$$
if $3\le k\le n_1$.
We see that $\psi_{u,u'}$ is the transposition of $u$ and $u'$.
Thus, \thetag{2-3} means that the transposition of $u$ and $u'$ {\it respects the tree structure}.
Given $p$ points $u^{(1)},\dots, u^{(p)}$, the permutation $\psi_{u^{(1)},u^{(2)}}\psi_{u^{(2)},u^{(3)}}\cdots\psi_{u^{(p-1)},u^{(p)}}$ is the  cycle $(u^{(1)},\dots, u^{(p)})$.
Hence, we can realize each cycle  $(u^{(1)},\dots, u^{(p)})$ in such a way that the cycle respects the tree structure.
In view of \thetag{2-2}, there is a partition $\Cal P$ of $X_2^{n_1}$ into atoms consisting of $p$ elements so that each atom of $\Cal P$ contains no more than one element of $Y_1$.
Every atom determines  a $p$-cycle.
Let $\gamma_1$ be the product of all these $p$-cycles.
Then $\gamma_1$ respects the tree structure, the order of $\gamma_1$ is $p$ and
the subsets
$$
\gamma_1^jY_1,\quad j=0,\dots,p-1,
$$
 of $X_2^{n_1}$ are mutually disjoint.
 Since $\gamma_1$ respects the tree structure, $g_1\in\Cal A$ (see (vii) for the definition of $g_1$).
 Pick an element $u_q$ from each atom $q\in\Cal P$.
 Denote by $\theta_1$ the product of  $p$ cycles generated by $(\gamma_1^ju_q)_{q\in\Cal P}$, $j=0,\dots,p-1$.
 Then  $\theta_1$  respects the tree structure, $\theta_1\gamma_1=\gamma_1\theta_1$ and
 the group generated by $\theta_1$ and $\gamma_1$ is transitive on $X_2^{n_1}$.

 %We now define a transformation $g_1$ of $X$ by setting
 %$$
% g_1(x_1,x_2,\dots,):=(\theta_0x_1,\gamma_1(x_2,\dots,x_{n_1}),x_{n_1+1},\dots).
% $$
 %It follows that $g_1$ is of order $\# X_1$ and the  the sets $\bigsqcup_{0\le j<\# X_1}g_1^j[x_1,Y_1]_1^{n_1}$, $x_1\in X_1$, are pairwise disjoint.

 On the second step,
by the Shannon-McMillan theorem, there is $n_2>n_1$  and a subset $Y_2\subset X_{n_1+1}^{n_2}$ such that
$$
\mu_{n_1+1}^{n_2}(Y_2)>1-\epsilon_2\quad\text{and}\quad \# Y_2\cdot p^{n_1-1}<p^{n_2-n_1}.
$$
 Therefore, there exists a bijection $\gamma_2\in\text{Bij}(X_{n_1+1}^{n_2})$ of order $p^{n_1-1}$ such that
the subsets
$$
\gamma_2^jY_2,\quad j=0,\dots,p^{n_1-1}-1,
$$
 of $X_{n_1+1}^{n_2}$ are mutually disjoint and $\gamma_2$ respects the tree structure.
 Then we construct a permutation $\theta_2$ of $X_{n_1+1}^{n_2}$ which commutes with
 $\gamma_2$ and jointly with $\gamma_2$ generate a transitive group of transformations
 of $X_{n_1+1}^{n_2}$.
 %We now define a transformation $g_2$ of $X$ by setting
% $$
% g_2(x_1,x_2,\dots,):=(x_1,\theta_1(x_2,\dots, x_{n_1}),\gamma_2(x_{n_1+1},\dots, x_{n_2}),x_{n_2+1},\dots).
% $$
%Then $g_2$ is of order $(\# X_1)^{n_1-1}$ and the subsets
%$\bigsqcup_{0\le j<\# X_1}g_2^j[y,Y_2]_1^{n_1}$, $y\in X_1^{n_1}$, are pairwise disjoint.

We continue this construction process infinitely many times to determine the entire sequence
$(n_k,Y_k,\theta_k, \gamma_k)_{k=1}^\infty$.
We obtain that
$$
\mu(X_1\times Y_1\times Y_2\times\cdots)=\mu_1^{n_1}(Y_1)\mu_{n_1+1}^{n_2}(Y_2)\cdots
\ge(1-\epsilon_1)(1-\epsilon_2)\cdots>0,
$$
as desired.
\endexample

\endcomment

\comment

Given $\delta>0$ and two subsets $A,B$ of $X$, we say that $A$ is $\delta$-full of $B$
if $\mu(A\cap B)>\delta\mu(A)$.

\proclaim{Theorem 2.5????} Let $\mu$ be the infinite product measure as above.
Let $G$ be minimal and $\mu$-finitary.
Then the dynamical system $(X,\mu,G)$ is ergodic.
\endproclaim
\demo{Proof}
Suppose that there is a $G$-invariant subset $A\subset X$ with $0<\mu(A)<1$.
Then there is $n>0$ and $x,y\in X_1^n$ such that the $n$-cylinder $[x]_1^n$ is $0.9999$-full of $A$ and the $n$-cylinder $[y]_1^n$ is $0.9999$-full of $A^c$.
Since $G$ is minimal, there is $g\in G$ with $g[x]_1^n=[y]_1^n$.
Since $G$ is $\mu$-finitary and $g$ is $\mu$-nonsingular, there exists $m>n$, a subset $Z\subset X_{n+1}\times\cdots\times X_m$ and a one-to-one mapping $Z\ni z\mapsto z'\in X_{n+1}\times\cdots\times X_m$ such that
\roster
\item"---" $[x]_1^n$ is 0.9999-full of $\bigsqcup_{z\in Z}[xz]_1^m$,
\item"---" $g[xz]_1^n=[yz']$ for each $z\in Z$,
\item"---" $[y]_1^n$ is 0.9999-full of $\bigsqcup_{z\in Z}[xz']_1^m$ and
\item"---"
$
g(x,z,x_{m+1}, x_{m+2},\dots)=(y, z', x_{m+1},x_{m+2}\dots)
$
if $z\in Z$.
\endroster
It follows from the latter equality and \thetag{2-1} that $\frac{d\mu\circ g}{d\mu}$ is constant on the cylinder $[xz]_1^m$ for each
$z\in Z$.
Let
$$
\align
E&:=\{v\in X_{n+1}^m\mid [xv]_1^m\text{ is 0.9-full of $A$}\}\text{ and }\\
F&:=\{v\in X_{n+1}^m\mid [yv]_1^m\text{ is 0.9-full of $A^c$}\}.
\endalign
$$
By a discrete version of Fubini theorem, $\mu_{n+1}^m(E)>0.9$
and $\mu_{n+1}^m(F)>0.9$.
Hence there is $z\in  E$ such that $z'\in F$.
It follows that $g[xz]_1^m=[yz']_1^m$.
Since $\frac{d\mu\circ g}{d\mu}$ is constant on $[xz]_1^m$ and
$[xz]_1^m$ is 0.9-full of $A$ and $A$ is invariant under $g$, it follows that $[yz']_1^m$
is 0.9-full of $A$, a contradiction.
\qed
\enddemo

\endcomment

\head 5. Ergodic measures for minimal transformation groups that are purely finitary with respect to the Haar measure
\endhead

\subhead 5.1. An ergodicity criterion for nonsingular transformation groups
\endsubhead
 In this subsection we provide a convenient sufficient condition for ergodicity of transformation groups on $X$.

\definition{Definition 5.1}  Let $\mu$ be a probability measure on $X$.
Let $G$ be a group of $\mu$-nonsingular transformations of $X$.
We say that $G$ is {\it compatible with $\mu$} if
for each %$\epsilon>0$ and
$N\in\Bbb N$, there exist $n>N$,
subsets $A_{y,y'}\subset [y]_1^n$  and elements $g_{y,y'}\in G$, $y,y'\in X_1^n$,
such that
$$
\align
\mu(A_{y,y'}) &>0.9\mu([y]_1^n),\\
g_{y,y'}A_{y,y'} &\subset[y']_1^n\quad\text{and}\\
\frac{d\mu\circ g_{y,y'}}{d\mu}(x) &=\frac{\mu([y']_1^n)}{\mu([y]_1^n)}\quad \text{at every $x\in A_{y,y'}$}
\endalign
$$
for all $y,y'\in X_1^n$.
\enddefinition

Of course, if $G$ is compatible with $\mu$ then $\mu(O)>0$ for each nonempty open subset $O\subset X$.

\proclaim{Proposition 5.2}
If $G$ is compatible with $\mu$ then $G$ is ergodic.
\endproclaim
\demo{Proof}
Let $F\subset X$ be a $G$-invariant Borel subset.
Suppose that $0<\mu(F)<1$.
Then we can find $n>0$, $y,y'\in X_1^n$, subset $A\subset [y]_1^n$ %and $B\subset [y']_1^n$
and an element $g\in G$ such that
$$
\align
 \mu(F\cap[y]_1^n) & >0.9\mu([y]_1^n),\\
  \mu(F^c\cap[y']_1^n) &>0.9\mu([y']_1^n),\\
 \mu(A)& >0.9\mu([y]_1^n),\\
 gA  & \subset[y']_1^n\quad\text{and}\\
\frac{d\mu\circ g}{d\mu}(x)& =\frac{\mu([y']_1^n)}{\mu([y]_1^n)}
\quad\text{at every $x\in A$.}
\endalign
$$
Then $\mu(F\cap A)>0.8\mu([y]_1^n)$, $g(F\cap A)\subset[y']_1^n$ and
$$
\mu(g(F\cap A))=\frac{\mu([y']_1^n)}{\mu([y]_1^n)}\mu(F\cap A)>0.8\mu([y']_1^n).
$$
Hence, $0=\mu(gF\cap F^c)\ge \mu(g(F\cap A)\cap F^c)>0.7\mu([y']_1^n)>0$,
a contradiction.
\qed
\enddemo

It follows from Proposition~5.2 and Claim~1.6 that if $G$ is compatible with $\mu$ then $G$ is minimal.

\subhead 5.2.  Uncountable collections of measures compatible with minimal groups of transformations
\endsubhead
We now show that  each minimal group of transformations which are purely finitary with respect to the Haar measure on $X$ is compatible with  a huge family of infinite product measures on $X$.

Let $\lambda_n$ stand for the equidistribution of $X_n$, $n\in\Bbb N$,
and let  $\lambda$ stand for the Haar measure on $X$.

\proclaim{Theorem 5.3} Let $G\subset \Cal A_\lambda^\bullet$ be a minimal group.
For each $n\in\Bbb N$, let $\beta_{n,0}$ and $\beta_{n,1}$ be two non-degenerated probabilities on $X_n$
such that for each $\omega\in\{0,1\}^\Bbb N$, the infinite product measure
$\bigotimes_{n\in\Bbb N}\beta_{n,\omega(n)}$ on $X$ is nonatomic.
Then there is an infinite subset $\Cal N\subset\Bbb N$ such that for each
$\omega\in\{0,1\}^\Bbb N$,
 $G$ is compatible with
 the infinite product measure $\lambda^\omega:=\bigotimes_{n\in\Bbb N}\gamma_n$ on $X$,
where
$$
\gamma_n:=
\cases
\lambda_n &\text{if $n\not\in\Cal N$},\\
\beta_{n,\omega(n)} & \text{if $n\in\Cal N$}.
\endcases
$$
Moreover, $G\subset\bigcap_{\omega\in\{0,1\}^\Bbb N}\Cal A^\bullet_{\lambda^\omega}$.
\endproclaim
\demo{Proof}
Enumerate the elements of $G$ as $g^{(1)},g^{(2)},\dots $.
We will construct $\Cal N$ inductively.
%Fix a decreasing sequence $(\epsilon_n)_{n=1}^\infty$ of positive reals with $\lim_{n\to\infty}\epsilon_n=0$.

The first step: since $G$ is minimal, there exist transformations $g_{y,y'}\in G$, $y,y'\in X_1$, such that
$g_{y,y'}[y]^1_1=[y']^1_1$.
Since $g^{(1)}\in\Cal A_\lambda^\bullet$ and $g_{y,y'}\in\Cal A_\lambda^\bullet$, there is $n_1>0$ and clopen subsets $A_{y,y'}\subset[y]_1^1$, $y,y'\in X_1$, such that
$$
\gather
\lambda(A_{y,y'})>0.9\lambda([y]_1^1),\\
 %g_{y,y'}A_{y,y'}\subset [y']_1,\\
 A_{y,y'}\text{ is the union of several $(n_1-1)$-cylinders and}\\
\text{$g^{(1)}$ and $g_{y,y'}$ change no more than $n_1-1$ first coordinates of each $x\in A_{y,y'}$,}
 \endgather
 $$
for all $y,y'\in X_1$.
Of course, $g^{(1)}$ and $g_{y,y'}$ preserve $\lambda$.

On the second step, for a given $\eta\in\{0,1\}$, replace $\lambda_{n_1}$ in the infinite product
$\bigotimes_{k=1}^\infty\lambda_k$ with $\beta_{n_1,\eta}$.
Then we obtain a new infinite product measure on $X$.
Denote it by $\lambda^\eta$.
Of course, $\lambda^\eta\sim\lambda$ for each $\eta\in\{0,1\}$.
%Hence $\Cal A_\lambda=\Cal A_{_0\lambda}=\Cal A_{_1\lambda}$.
Since $G$ is minimal, there exist transformations $g_{y,y'}\in G$, $y,y'\in X_1^{n_1}$, such that
$g_{y,y'}[y]_1^{n_1}=[y']_1^{n_1}$.
Since $g^{(1)}$, $g^{(2)}$ and $g_{y,y'}$ belong to $\Cal A_\lambda^\bullet=\Cal A_{\lambda^\eta}^\bullet$ for each $\eta\in\{0,1\}$, there is $n_2>n_1$ and clopen subsets $A_{y,y'}\subset[y]_1^{n_1}$, $y,y'\in X_1^{n_1}$, such that
$$
\gather
\min_{\eta=0,1} {\lambda^\eta}(A_{y,y'})>0.99\lambda^\eta([y]_1^{n_1}),\\
% g_{y,y'}A_y\subset [y']_1^{n_1},\\
 A_{y,y'}\text{ is the union of several $(n_2-1)$-cylinders,}\\
\text{$g^{(1)}, g^{(2)}$ and $g_{y,y'}$ change no more than $n_2-1$ first coordinates of each $x\in A_{y,y'}$}.\\
 \endgather
 $$
 It follows that for all $y,y'\in X_1^{n_1}$,
 $$
\frac{d\lambda^\eta\circ g_{y,y'}} {d\lambda^\eta}(x)=
\frac{{\lambda^\eta}([y']_1^{n_2-1})}{{\lambda^\eta}([y]_1^{n_2-1})}=
\frac{{\lambda^\eta}([y']_1^{n_1})}{{\lambda^\eta}([y]_1^{n_1})}\quad\text{at  each $x\in A_{y,y'}$}
$$
because $\lambda^\eta$ is the infinite product measure whose $j$-th coordinate is the equidistribution on $X_j$
whenever $n_1<j<n_2$.

On the $l$-th step, for a given $\eta\in\{0,1\}^{l-1}$, replace $\lambda_{n_1},\dots,\lambda_{n_{l-1}}$ in the product
$\bigotimes_{n=1}^\infty\lambda_n$ with $\beta_{n_1,\eta(1)}$,\dots, $\beta_{n_{l-1},\eta(l-1)}$ respectively.
We then obtain a new infinite product measure on $X$.
Denote it by $\lambda^\eta$.
Of course, $\lambda^\eta\sim\lambda$ for each $\eta\in\{0,1\}^{l-1}$.
%Hence $\Cal A_\lambda=\Cal A_{_0\lambda}=\Cal A_{_1\lambda}$.
Since $G$ is minimal, there exist transformations $g_{y,y'}\in G$, $y,y'\in X_1^{n_{l-1}}$, such that
$g_{y,y'}[y]_1^{n_{l-1}}=[y']_1^{n_{l-1}}$.
Since $g^{(j)},g_{y,y'}\in\Cal A_\lambda^\bullet=\Cal A_{\lambda^\eta}^\bullet$ for each $j=1,\dots,l$ and $\eta\in\{0,1\}^{l-1}$, there is $n_l>n_{l-1}$ and clopen subsets $A_{y,y'}\subset[y]_1^{n_{l-1}}$, $y,y'\in X_1^{n_{l-1}}$, such that
$$
\gather
\min_{\eta \in\{0,1\}^{l-1}} \lambda^\eta(A_{y,y'})>(1-10^{-l})\lambda^\eta([y]_1^{n_l-1}),
\tag5-1\\
   A_{y,y'}\text{ is the union of several $(n_l-1)$-cylinders,}
  \tag5-2\\
% \text{ if $x=(x_j)_{j=1}^\infty\in A_y$ then $g_{y,y'}x=(\dots, x_{n+1},x_{n+2},\dots)$.}
\text{$g^{(j)},g_{y,y'}$ change no more than $n_l-1$ first coordinates of each $x\in A_{y,y'}$}\tag5-3\\
 \endgather
 $$
for all $j=1,\dots,l$ and $y,y'\in X_1^{n_{l-1}}$.
Then we have that
$$
\gather
g_{y,y'}A_{y,y'}\subset [y']_1^{n_l-1},\tag5-4\\
\text{and}\quad\frac{d\lambda^\eta\circ g_{y,y'}} {d\lambda^\eta}(x)=
\frac{\lambda^\eta([y']_1^{n_{l-1}})}{\lambda^\eta([y]_1^{n_{l-1}})}\quad\text{at  each $x\in A_{y,y'}$}\tag5-5
\endgather
$$
We now let $\Cal N:=\{n_1,n_2,\dots\}$.
Take an arbitary $\omega\in\{0,1\}^\Bbb N$.
Let $\lambda^\omega$ be the measure defined in the statement of the theorem.

We first claim that $G\subset \Cal A^\bullet_{\lambda^\omega}$.
Indeed, for each $l>0$, let
$$
B_l:=\bigsqcup_{y\in X_1^{n_l-1}}A_{y,y}.
$$
Denote by $\omega| (l-1)$ the restriction of $\omega$ to $\{1,2,\dots,l-1\}$.
Thus, $\omega|(l-1)\in\{0,1\}^{l-1}$.
Then in view of \thetag{5-1} and \thetag{5-2},
$$
\multline
\lambda^\omega(B_l)=\lambda^{\omega\restriction(l-1)}(B_l)=
\sum_{y\in X_1^{n_l-1}}\lambda^{\omega\restriction(l-1)}(A_{y,y})>
(1-10^{-l})\sum_{y\in X_1^{n_l-1}}\lambda^\omega([y]_1^{n_l-1})\\
=1-10^{-l}.
\endmultline
$$
From this, \thetag{5-2} and \thetag{5-3} we deduce that
$g^{(1)},g^{(2)},\dots\in \Cal A^\bullet_{\lambda^\omega}$, as desired.
In particular, $G$ is $\lambda^\omega$-nonsingular.

It remains to
show
 that $G$ is compatible with $\lambda^\omega$.
To this end we fix  $l\in\Bbb N$.
Let $(A_{y,y'})_{y,y'\in X_1^{n_{l-1}}}$ and $(g_{y,y'})_{y,y'\in X_1^{n_{l-1}}}$ be  the objects that were
 determined on the $l$-th step.
 It follows from \thetag{5-2} and \thetag{5-1} that
 $$
 \lambda^\omega(A_{y,y'})=\lambda^{\omega |(l-1)}(A_{y,y'})>0.9\lambda^{\omega |(l-1)}([y]_1^{n_{l-1}})
 =0.9\lambda^\omega([y]_1^{n_{l-1}})
 \tag5-6
 $$
 for all $y,y'\in X_1^{n_{l-1}}$.
Since \thetag{5-3} and \thetag{5-5} hold for $g_{y,y'}$, it follows  that
$$
\frac{d\lambda^\omega\circ g_{y,y'}} {d\lambda^\omega}(x)
=\frac{d\lambda^{\omega |(l-1)}\circ g_{y,y'}} {d\lambda^{\omega |(l-1)}}(x)
=
\frac{\lambda^{\omega |(l-1)}([y']_1^{n_{l-1}})}
{\lambda^{\omega |(l-1)}([y]_1^{n_{l-1})}}
=\frac{\lambda^{\omega }([y']_1^{n_{l-1}})}{\lambda^{\omega }([y]_1^{n_{l-1}})}
$$
for all $y,y'\in X_1^{n_{l-1}}$.
This equality, \thetag{5-2}--\thetag{5-4} and \thetag{5-6} yield  that $G$ is compatible with $\lambda^\omega$.
\qed
\enddemo

\subhead 5.3. Kakutani theorem and a corollary from Theorem~5.3
\endsubhead
We remind that given two probability measures $\alpha,\beta$ on a finite set $B$,
the {\it Hellinger distance} $H(\alpha,\beta)$ between $\alpha$ and $\beta$ is
$$
H(\alpha,\beta):=\frac1{\sqrt 2}\sqrt{\sum_{b\in B}\Big(\sqrt{\alpha(b)}-\sqrt{\beta(b)}\Big)^2}.
$$
It is easy to verify that
$$
1-H^2(\alpha,\beta)=\sum_{b\in B}\sqrt{\alpha(b)\beta(b)}.
$$
We remind the well know dichotomy result on  infinite product measures.

\proclaim{Kakutani theorem \cite{Ka}} For each $n\in\Bbb N$, fix  two nondegenerated probabilities $\mu_n$ and $\nu_n$ on $X_n$.
Let $\mu:=\bigotimes_{n=1}^\infty\mu_n$ and
 $\nu:=\bigotimes_{n=1}^\infty\nu_n$.
We have that
 \roster
 \item"---"
 if $\prod_{n=1}^\infty(1-H^2(\mu_n,\nu_n))>0$
 %\sum_{x_n\in X_n}\sqrt{\mu_n(x_n)\lambda_n(x_n)}\ne 0$
  then $\mu\sim\nu$,
  \item"---"
 if $\prod_{n=1}^\infty(1-H^2(\mu_n,\nu_n))= 0$ then $\mu\perp\nu$.
 \endroster
\endproclaim

Given two points $\omega,\omega'\in\{0,1\}^\Bbb N$,
we write $\omega\sim\omega'$ if they are tail  equivalent (see Remark~3.3 for the definition).
We obtain the following corollary from the Kakutani theorem, Theorems~5.3 and Proposition~5.2.

\proclaim{Corollary 5.4} Let $G\subset\Cal A_\lambda^\bullet$ be a minimal group and let $\delta>0$.
Let $\beta_{n,0}$ and $\beta_{n,1}$ be two non-degenerated probabilities on $X_n$
such that $H(\beta_{n,0},\beta_{n,1})>\delta$,  $H(\beta_{n,0},\lambda_n)>\delta$, $H(\beta_{n,1},\lambda_n)>\delta$ and for each $\omega\in\{0,1\}^\Bbb N$, the infinite product measure
$\bigotimes_{n\in\Bbb N}\beta_{n,\omega(n)}$ on $X$ is nonatomic.
Given $\omega\in\{0,1\}^\Bbb N$, let $\lambda^\omega$ stand for the infinite product measure defined in the statement
of Theorem~5.3.
Then the following are satisfied:
\roster
\item"(i)" $G$ is $\lambda^\omega$-nonsingular and ergodic for each
$\omega\in\{0,1\}^\Bbb N$.
\item"(ii)"
If $\omega,\omega'\in\{0,1\}^\Bbb N$ and $\omega\not\sim\omega'$
then $\lambda^\omega\perp \lambda^{\omega'}$.
\item"(iii)"
If  $\omega,\omega'\in\{0,1\}^\Bbb N$ and $\omega\sim\omega'$
then $\lambda^\omega\sim \lambda^{\omega'}$.
\item"(iv)"  $\lambda^\omega\perp \lambda$ for each $\omega\in\{0,1\}^\Bbb N$.
\endroster
\endproclaim

\demo{Proof} Proposition 5.2 and Theorem 5.3 imply (i).
We note that for each pair $\omega,\omega'\in\{0,1\}^\Bbb N$,
$$
\sum_{n=1}^\infty H^2((\lambda^\omega)_n,(\lambda^{\omega'})_n)=
\sum_{n\in\Cal N} H^2((\lambda^\omega)_n,(\lambda^{\omega'})_n)=
\sum_{n\in\Cal N} H^2(\beta_{n,\omega(n)},\beta_{n,\omega'(n)}).
$$
Hence, the Kakutani theorem implies
 (ii) and (iii).
 In a similar one, one can prove~(iv).
\qed
\enddemo

\comment

\proclaim{Theorem 2.10} Let $G\subset\Cal A_0$ be a minimal group.
Then there is an infinite subset $\Cal N\subset\Bbb N$ such that for each sequence $(\mu_n)_{n=1}^\infty$ with $\mu_n$ being the equidistribution on $X_n$ for each $n\not\in\Cal N$, the infinite product $\mu:=\bigotimes_{n\in\Bbb N}\mu_n$ is compatible with $G$.
\endproclaim

\demo{Proof} We will construct $\Cal N$ inductively.
Let $n_1:=1$.
Since $G$ is minimal, there is a mapping $g_1:X_1\times X_1\to G$ such that
$g_1(y,y')[y]_1=[y']_1$ for all $y,y'\in X_1$.
Since $g_1(y,y')\in\Cal A_0$, there is $n_2>n_1$ such that   $g_1(y,y')$ does not change any coordinate greater than $n_2-1$ for all $y,y'\in X_1$.
On the second step we select
a mapping $g_2:X_1^{n_2}\times X_1^{n_2}\to G$ such that
$g_2(y,y')[y]_1=[y']_1$ for all $y,y'\in X_1^{n_2}$.
Then there is $n_3>n_2$ such that
$g_2(y,y')$ does not change any coordinate greater than $n_3-1$ for all $y,y'\in X_1^{n_2}$.
We continue this process infinitely many times to construct the entire sequence $(n_k)_{k=1}^\infty$.
Let $\Cal N:=\{n_k\mid k\in\Bbb N\}$.
It is straightforward to verify that $\Cal N$ is as desired.
\qed
\enddemo

\endcomment

\head 6. Irreducible Koopman representations for
branch groups
\endhead

\comment

Let $\mu$ be a  $G$-quasiinvariant measure  on $X$.
 We remind that the corresponding
Koopman representation of $G$ in $L^2(X,\mu)$
is denoted by $\kappa_\mu=(\kappa_\mu(g))_{g\in G}$

\endcomment

The following is a well-known result (see e.g. \cite{Ta, Corollary 16}).

\proclaim{Theorem on irreducibility of Koopman representations} Let $H$ be an ergodic group of nonsingular transformations of a probability space $(Y,\nu)$.
Denote by $\kappa_\nu$ the corresponding Koopman representation of $H$ in $L^2(Y,\nu)$.
Given $f\in L^\infty(Y,\nu)$, let $L_f$ denote the linear operator (in the Hilbert space $L^2(Y,\nu)$) of multiplication by $f$.
 If, for each $f\in L^\infty(Y,\nu)$, the operator $L_f $ belongs to the von Neumann algebra generated by the Koopman operators $\kappa_\nu(h)$, $h\in H$,  then $\kappa_\nu$ is irreducible.
\endproclaim

We now return to the notation of the previous section.
Fix a subgroup $G\subset\Cal A$.
For a cylinder $O$  in $X$,  let
$$
\align
%G_O^+&:=\{g\in G\mid g O=O\},\\
G_O&:=\{g\in G\mid \text{Supp\,}g\subset O\}\quad\text{and}\\
\Cal I_O&:=\{h\in L^2(X,\mu)\mid \kappa_\mu(g)h=h\text{ for each $g\in G_O$}\},
\endalign
$$
where Supp\,$g:=\{x\in X\mid gx\ne x\}$.
%The group $G_O$ is called {\it the rigid stabilizer of $O$}.

\comment

\proclaim{Proposition 6.1} Let $g\in G$.
\roster
\item"(i)"
If $x\in O$ and $gx\in O$ then $g\in G_O^+$.
\item"(ii)"
If $x\in O$ then $(Gx)\cap O=G_O^+x$,
\item"(iii)"
If $G$ is minimal then $G_O^+$ is minimal on $O$.
\endroster
\endproclaim

\demo{Proof}
We note that $g$ has a special structure \thetag{1-1}  because  $g\in \Cal A$.
Therefore, a routine verification shows that (i) holds.
Of course, (i) implies (ii) directly.
In turn, (iii)  follows from (ii) straightforwardly.
\qed
\enddemo

\endcomment

 As above, $\lambda$ denotes the Haar measure on $X$.
We remind that $G$ is called branch if it is minimal and for each $n>0$,
the subgroup
$
\prod_{u\in X_1^n}G_{[u]_1^n}
$
is of finite index in $G$.
%Experts in the field of rooted trees denote this subgroup by rist$_G(n)$.

\comment
We introduce a concept of quasi-branch.

\definition{Definition 6.2} A countable group $G\subset\Cal A$ is called {\it quasi-branch} if $G$ is minimal and the group
$
\prod_{u\in X_1^n}G_{[u]_1^n}
$
is nontrivial and of finite index in $G_{[u]_1^n}^+$
for each $u\in X_1^n$ and $n\in\Bbb N$.
\enddefinition

Of course, a branch group is quasi-branch.
A quasi-branch is weakly branch.

\endcomment

The following lemma is an analogue of \cite{DuGr1, Proposition~5}, where it was stated under different conditions. It was assumed there that  the group  $G$ is weakly branch, the tree is regular, i.e. $X_1=X_2=\cdots$ and
$\lambda$ is a Bernoulli measure on $X$.
The proofs \cite{DuGr1, Proposition~5} and Lemma~6.1 below are totally different.

We will use the following notation: for a subset $A$ of a set $B$, by $A^c$ we denote the complement of $A$ in $B$.

\proclaim{Lemma 6.1}  Let $\mu$ be a nonatomic probability on $X$
%$\mu=\bigotimes_{n=1}^\infty\mu_n$
and let $\mu\perp\lambda$.
Let $G$ be  branch and $\mu$-nonsingular.
Then for each cylinder $O$ in $X$,
$$
\Cal I_O=\{h\in L^2(X,\mu)\mid \text{\rom{Supp}\,}h\subset O^c\}.
$$
\endproclaim
\demo{Proof}
%Since $G$ is $\mu$-finitary, $G$ is $\mu$-nonsingular.
Take $h\in\Cal I_O$.
Then for each $g\in G_O$,
$$
\big(\kappa_\mu(g)h\big)(x)=h(g^{-1}x)\sqrt{\frac{d\mu\circ g^{-1}}{d\mu}(x)}=h(x)
$$
at a.e. $x\in X$.
Let $\text{Supp}\, h:=\{x\in X\mid h(x)\ne 0\}$.
It follows that $\text{Supp}\,h$ is invariant under $G_O$ and
$$
\frac{d\mu\circ g^{-1}}{d\mu}(x)=\frac{|h(x)|^2}{|h(g^{-1}x)|^2}\quad\text{for each $g\in G_O$ at a.e.
$x\in \text{Supp}\,h.$}
$$
Hence $ |h^2|\cdot\mu$ is a $G_O$-invariant finite measure\footnote{By $|h^2|\cdot\mu$ we mean
the measure on $X$ which is absolutely continuous with respect to $\mu$ and such that $\frac{d(|h^2|\cdot\mu)}{d\mu}=|h^2|$.}
 which is absolutely continuous with respect to $\mu$.
This implies, in turn, that $1_{O\cap \text{Supp}\,h}\cdot |h^2|\cdot\mu$ is a $G_O$-invariant finite measure which is absolutely continuous with respect to $1_O\cdot\mu$.
Since $O$ is a cylinder, there
  exist some $n>0$ and $v\in X_1^n$ such that
$O=[v]_1^n$.
We let $\widetilde G_O:=\prod_{u\in X_1^n}G_{[u]_1^n}$.
It follows (almost tautologically)  that $1_{O\cap \text{Supp}\,h}\cdot |h|^2\cdot\mu$ is  $\widetilde G_O$-invariant.
Since $G$ is branch,  %$\widetilde G_O$ is a subgroup of finite index in $G$.
%Since $\widetilde G_O\subset G_O^+\subset G$,
%it follows that
$\widetilde G_O$ is a subgroup of finite index in $G$. %$G_O^+$.
Hence, there are elements $g_1,\dots,g_l\in G$ %G_O^+$
such that $G=\bigsqcup_{j=1}^l\widetilde G_Og_j$.
The sing $\ll$ below is used to denote  the absolute continuity of measures.
As for each $j$,
$$
(1_{O\cap \text{Supp}\,h}\cdot |h|^2\cdot\mu)\circ g_j\ll \mu\circ g_j
\sim \mu,
$$
we obtain that
$$
\sum_{j=1}^n(1_{O\cap \text{Supp}\,h}\cdot |h|^2\cdot\mu)\circ g_j\ll\mu.
\tag6-1
$$
The left-hand side measure in this relation is finite.
It is straightforward to verify that it is  invariant under $G$. %$G_O^+$.
Since $G$ in minimal on $X$, it follows from Claim~1.5 that
$G$ is uniquely ergodic.
%We can identify naturally $O$ with the infinite product space $X_{n+1}\times X_{n+2}\times\cdots$.
%Then $O$ is $G_O^+$-uniquely ergodic in view of Claim~1.5.
Hence, each nontrivial finite $G$-invariant measure  is proportional to $\lambda$.
Therefore, it follows from \thetag{6-1} that $\lambda\ll\mu$.
This contradicts to the assumption of the lemma.
Hence,  the left-hand side of \thetag{6-1} is trivial.
This implies, in turn, that  $\mu(O\cap \text{Supp}\,h)=0$, i.e. Supp$\,h\subset O^c$.
Thus,
$$
\Cal I_O\subset\{h\in L^2(X,\mu)\mid \text{\rom{Supp}\,}h\subset O^c\}.
$$
The converse inclusion is obvious.
\qed
\enddemo

\proclaim{Theorem 6.2} Let $\mu$ be a nonatomic  measure on $X$.
%$\mu:=\bigotimes_{n=1}^\infty\mu_n$.
   Let $G$ be  branch,  $\mu$-nonsingular and ergodic.
If  $\mu\perp\lambda$
then
 $\kappa_\mu$
is irreducible.
\comment

Let $\lambda$ be another (infinite product) measure $\lambda$ such that $G$ is $\lambda$-finitary.
Denote by $\kappa_\mu$ and $\kappa_\lambda$ the Koopman represenatations
generated by the nonsingular systems $(X,\mu, G)$  and $(X,\lambda , G)$.
Then
\roster
\item"---"
If $\mu\perp\lambda$ then $\kappa_\mu$ and $\kappa_\lambda$ are unitarily nonequivalent,
\item"---"
If $\mu\sim\lambda$ then $\kappa_\mu$ and $\kappa_\lambda$ are unitarily equivalent.
\endroster

\endcomment
\endproclaim

\demo{Proof}
Denote by $\Cal N_\mu$ the von Neumann algebra generated by $\kappa_\mu$,
i.e.
$$
\Cal N_\mu=\{\kappa_\mu(g)\mid g\in G\}''.
$$
Fix a cylinder $O\subset X$.
Then the orthogonal projection $P$ onto $\Cal I_O$ belongs to $\Cal N_\mu$.
It follows from Lemma~6.1 that $P=L_{1_{O^c}}$.
Hence, $L_{1_{O^c}}\in \Cal N_\mu$.
This yields  that $L_{1_{O}}=I-L_{1_{O^c}}\in \Cal N_\mu$.
Therefore, $L_{1_{U}}\in \Cal N_\mu$ for each clopen subset $U\in X$.
This, in turn, implies that $L_f\in \Cal N_\mu$ for each $f\in L^\infty(X,\mu)$.
It remains to apply the theorem on irreducibility of Koopman representations.
\qed
\enddemo

A natural question related to Theorem~6.2 arises: whether $\kappa_\lambda\restriction L^2_0(X,\lambda)$ is irreducible?
In the case of regular trees,
the following claim was proved in \cite{BaGr, Theorem 3}.
We extend it to the general case.

\proclaim{Proposition 6.3} Let $G$ be an arbitrary minimal countable subgroup of $\Cal A$.
Then $\kappa_\lambda$ splits into  a countable orthogonal sum of finitely dimensional representations of $G$.
\endproclaim
\demo{Proof}
Let $K$ denotes the closure of $G$ in $\Cal A$.
Since $\Cal A$ is compact, $K$ is compact.
Of course, $K$ acts on $X$ and preserves $\lambda$.
Hence $\kappa_\lambda$ extends to the Koopman representation of $K$.
We denote the ``extended'' representation by the same symbol $\kappa_\lambda$.
Since $G$ is dense in $K$, it follows that $(\kappa_\lambda(k))_{k\in K}$ is irreducible
if and only if $(\kappa_\lambda(g))_{g\in G}$ is irreducible.
Since $G$ is minimal, $(X,\lambda, G)$ is ergodic.
Hence, the dynamical system $(X,\lambda, K)$ is ergodic.
Each ergodic action of a compact group is  transitive action (mod 0).
Thus, there is a closed subgroup $K_0\subset K$ such that $(X,\lambda, K)$
is isomorphic to $(K/K_0, \xi,K)$, where $K$ acts on the coset space $K/K_0$ by left translations
and $\xi$ is the Haar measure on $K/K_0$.
It follows from the Peter-Weyl theorem that the corresponding Koopman representation
splits into a countable orthogonal sum of finitely dimensional irreducible unitary represenations of $K$.
The restriction of each of these $K$-irreducible representations to $G$ is an irreducible unitary representation of $G$.
\qed
\enddemo

\comment

We note that since $\mu$ is finite,  each $G$-invariant subset $B$ in $X$
generates an $\kappa_\mu$ invariant closed subspace $1_B\cdot L^2(X,\mu)$
in $L^2(X,\mu)$.
Therefore, we obtain the following corollary from Theorem~3.2.

\proclaim{Corollary 3.3} Under the condition of Theorem~3.2, $G$ is ergodic.
\endproclaim

\endcomment

\proclaim{Theorem 6.4} Let $\mu$ and $\mu'$ be two
%infinite product
nonatomic measures on $X$ such that
$\mu\perp\lambda$ and $\mu'\perp\lambda$.
Let $G$ be  branch, nonsingular  and ergodic  with respect to both $\mu$ and  $\mu'$.
The following holds:
\roster
\item"(i)"
If $\mu\sim\mu'$ then $\kappa_\mu$ and $\kappa_{\mu'}$ are unitarily equivalent,
\item"(ii)"
If $\mu\perp\mu'$ then
$\kappa_\mu$ and $\kappa_{\mu'}$ are  disjoint (i.e. not unitarily equivalent).
\endroster
\endproclaim
\demo{Proof}
(i) is obvious.

(ii) Suppose that $\mu\perp\mu'$ but
$\kappa_\mu$ and $\kappa_{\mu'}$ are   unitarily equivalent.
Then there exists a unitary isomorphism $V:L^2(X,\mu)\to L^2(X,\mu')$  such that
$$
V\kappa_{\mu}(g)V^{*}=\kappa_{\mu'}(g)\quad\text{for each $g\in G$.}
$$
For a function $f\in L^\infty(X,\mu)$, we denote by $L_f$ the operator of multiplication by $f$ in $L^2(X,\mu)$.
In a similar way, for a function $f\in L^\infty(X,\mu')$, we denote by $L_f'$ the operator of multiplication by $f$ in $L^2(X,\mu')$.
Let
$$
\align
\Cal I_O&:=\{h\in L^2(X,\mu)\mid \kappa_\mu(g)h=h\text{ for all $g\in G_O$}\}\quad\text{and}\\
\Cal I_O'&:=\{h\in L^2(X,\mu')\mid \kappa_{\mu'}(g)h=h\text{ for all $g\in G_O$}\}
\endalign
$$
and let $P_O$ and $P_O'$ stand for the orthogonal projections onto $\Cal I_O$
and $ \Cal I_O$ respectively.
It is straightforward to verify that $V\Cal I_O=\Cal I_O'$.
Hence $VP_O=P_O'V$.
By Lemma~6.1, $P_O=L_{1_{O^c}}$ and $P_O'=L_{1_{O^c}}'$.
Hence, $VL_{1_{O^c}}V^*=L_{1_{O^c}}'$.
As $O$ is an arbitrary cylinder in $X$, we deduce that
$$
VL_fV^*=L_f'
$$
 for each bounded Borel function $f:X\to\Bbb C$.
However, this is only possible if $\mu\sim\mu'$ and $Vf=f\cdot\frac{d\mu'}{d\mu}$,
a contradiction.
\comment

$\Cal I_O:=\{h\in L^2(X,\mu)\mid\text{Supp\,}h\subset O^c\}$

Let $\Cal N_\mu$ and $\Cal N_{\mu'}$ be the von Neumann algebras generated by
$\kappa_\mu$ and $\kappa_{\mu'}$ respectively.
%It follows from  the proof of Theorem~6.4 that $L^\infty(X,\mu)\subset \Cal N_\mu$
and
$L^\infty(X,\mu')\subset \Cal N_{\mu'}$.
It follows from the proof of Theorem~6.5 that the ``conjugacy by $V$'', i.e. the mapping
$$
\kappa_{\mu}(g)\mapsto V\kappa_{\mu}(g)V^{-1},\quad g\in G,
$$
 extends uniquely to the $*$-isomorphism of $\Cal N_\mu$ onto $\Cal N_{\mu'}$.
Fix a cylinder $O$ in $X$.
Then $L_{1_{O^c}}$ is the minimal projector in $\Cal N_\mu$ that commutes with $\kappa_\mu(g)$ for each $g\in G_O$.
Hence $VL_{1_{O^c}}V^{-1}$ is the minimal projector in $\Cal N_{\mu'}$ that commutes with $V\kappa_{\mu}(g)V^{-1}=\kappa_{\mu'}(g)$ for each $g\in G_O$.
However, this projector is exactly $L_{1_{O^c}}$ in $L^2(X,\mu')$.
Thus,
$$
VL_{1_{O^c}}V^{-1}=L_{1_{O^c}}.\tag6-2
$$
Since $O$ is arbitrary, $V$ admits a spatial realization, i.e. there is $\phi:X\to X$ such that
$\mu\circ\phi^{-1}\sim\mu'$ and
$Vf:=f\circ\phi^{-1}\sqrt{\frac{d\mu'\circ\phi}{d\mu}}$.
Hence, \thetag{6-2} yields that $1_{O^c}\circ\phi=1_{O^c}$ for each cylinder $O$.
This is only possible if $\phi(x)=x$ for $\mu$-a.e. $x\in X$.
Hence, $\mu\sim\mu'$, a contradiction.
\endcomment
\qed
\enddemo

Denote by $\Cal P(X)$ the set of probability  measures on $X$.
Endow $\Cal P(X)$ with the  weak topology.
Then $\Cal P(X)$ is a compact Polish affine space.
Let $\Cal P_\infty(X)$  denote the subset of infinite product probabilities on $X$.
Then $\Cal P_\infty(X)$ is a closed subset of $\Cal P(X)$.

\proclaim{Theorem 6.5} Let $G$ be a branch subgroup of $\Cal A_\lambda^\bullet$.
Then there is a continuous mapping $\{0,1\}^\Bbb N\ni\omega\mapsto\mu_\omega\in\Cal P_\infty(X)$
such that:
\roster
\item"(i)" $G$ is nonsingular and ergodic with respect to $\mu_\omega$ for each $\omega\in \{0,1\}^\Bbb N$.
\item"(ii)" The Koopman representation $\kappa_{\mu_\omega}$ of $G$ is irreducible for each $\omega\in \{0,1\}^\Bbb N$.
\item"(iii)"
If $\omega,\omega'\in\{0,1\}^\Bbb N$ and $\omega\sim\omega'$ then
$\kappa_{\mu_\omega}$ and $\kappa_{\mu_\omega}$ are unitarily equivalent.
\item"(iv)"
If $\omega,\omega'\in\{0,1\}^\Bbb N$ and $\omega\not\sim\omega'$ then
$\kappa_{\mu_\omega}$ and $\kappa_{\mu_\omega}$ are disjoint (not unitarily equivalent).
\endroster
\endproclaim

\demo{Proof}
As $G$ is minimal and $G\subset\Cal A_\lambda^\bullet$, we apply Corollary~5.4 to define  the mapping
$\{0,1\}^\Bbb N\ni\omega\mapsto\mu_\omega\in\Cal P_\infty(X)$ satisfying the conditions (i)--(iii)
of that corollary.
It is routine to verify that this mapping is continuous.
Thus, (i) follows from Corollary~5.4.
Theorem~6.2 implies (ii).
The properties (iii) and (iv) follow from~Theorem~6.4.
\qed
\enddemo

\head 7.  Weak containment of unitary
representations
\endhead

Let $U=(U_g)_{g\in G}$ and $V=(V_g)_{g\in G}$ are unitary representations of $G$ in a separable Hilbert
spaces $\Cal H$ and $\Cal D$ respectively.
We remind that $U$ is {\it weakly contained in} $V$ (denoted by $U\prec V$) if for each finite subset $F\subset G$, vector $h\in \Cal H$ and $\epsilon>0$ there exist finitely many vectors $d_1,\dots,d_k\in D$ such that
$$
\max_{g\in F}\Bigg|\langle  U_gh,h\rangle-\sum_{j=1}^k\langle V_gd_j,d_j\rangle \Bigg|<\epsilon,
$$
where $\langle .,.\rangle$ means the inner product.
For instance, $G$ is amenable if and only if the trivial representation of $G$ is weakly contained in the left regular representation of $G$.
If  $U\prec V$ and $V\prec U$ then $U$ and $V$ are called {\it weakly equivalent}.

\proclaim{Theorem 7.1}  Given $n>0$, let $X_n$ be a finite set and let   $\mu_n$ and $\nu_n$ be  two non-degenerated probabilities  on $X_n$.
Let $X:=X_1\times X_2\times\cdots$, $\mu:=\mu_1\otimes\mu_2\otimes\cdots$,
$\nu:=\nu_1\otimes\nu_2\otimes\cdots$.
Suppose that $\mu$ and $\nu$ are nonatomic.
Let $G$ be a countable subgroup of $\Cal A_\mu^\bullet$.
%uch that $G\subset \Cal A_\mu^\bullet\cap\Cal A_\nu^\bullet$.
%Suppose that for each $g\in G$,
%$$
%\lim_{k\to\infty}\nu_1^n(Y_{g,k}^c)\max_{y\in X_1^k}\frac{\mu_1^k(y)}{\nu_1^k(y)}= 0.\tag7-1
%$$
Then    $\kappa_\mu\prec\kappa_\nu$.
Hence, if
 $G\subset \Cal A_\mu^\bullet\cap\Cal A_\nu^\bullet$ then
$\kappa_\mu$ and $\kappa_\nu$ are weakly equivalent.
\endproclaim

\demo{Proof} Fix a finite subset $F\subset G$,  $\epsilon>0$,
 $n
 \in\Bbb N$ and a function $f:X_1\times\cdots\times X_n\to\Bbb C$.
 By $f\otimes 1$ we denote the mapping from $X\to\Bbb C$ such that
 $$
 f\otimes 1(x_1,x_2,\dots)=f(x_1,\dots,x_n)
  $$
  for each $x=(x_l)_{l=1}^\infty\in X$.
 Since $F\subset\Cal A$,
 for each $g\in F$ and $m>0$, there exist
 $$
 g_m\in \text{Bij}(X_1^m)\quad\text{and
 mappings}\quad
 \alpha_{g,l}:X_1^l\to\text{Bij}(X_{l+1}),\quad l=m,m+1,\dots,
 $$
 such that
 $$
 gx=(g_m(x_1,\dots,x_m), \alpha_{g,m}(x_1,\dots,x_m)[x_{m+1}],\alpha_{g,m+1}(x_1,\dots,x_{m+1})[x_{m+2}],\dots).
$$
 For $k>0$, we let
 $$
 \align
Y_{g,k}&:=\{y\in X_1^k\mid \alpha_{g,l}(y,z)=I\text{ for all $z\in X_{k+1}^l$ and  each $l>k$}\}
\quad\text{and}\\
Y_{k}&:=\bigcap_{g\in F}Y_{g,k}\subset X_1^k.
\endalign
$$
Then $\lim_{k\to\infty}\mu_1^k(Y_{g,k})
 =1$ because
 $F\subset \Cal A_\mu^\bullet$.  %\cap\Cal A_\nu^\bullet$.
    Hence, there is $k_\epsilon>n$ such that  $\mu_1^{k_\epsilon}(Y_{k_\epsilon})>1-\epsilon$
    and $\mu_1^{k_\epsilon}(g^{-1}Y_{k_\epsilon})>1-\epsilon$ for each $g\in F$.
    Let $\widetilde Y:=\bigsqcup_{y\not\in Y_{k_\epsilon}}[y]_1^{k_\epsilon}$.
    Then for each $g\in F$,
    by the Cauchy-Schwarz inequality,
    $$
    \int_{\widetilde Y}(f\circ g_n^{-1}\cdot f)(x_1,\dots,x_n)\sqrt{\frac{d\mu\circ g^{-1}}{d\mu}}(x)\,d\mu(x)
    <C^2\sqrt{\mu(\widetilde Y)\mu(g^{-1}\widetilde Y)},
    $$
    where $x_1,\dots, x_n$ are the first $n$ coordinates of $x\in \widetilde Y$ and   $C:=\max_{z\in X_1^n}|f(z)|$.
    We have that $\mu(\widetilde Y)
    =\mu_1^{k_\epsilon}(Y_{k_\epsilon}^c)<\epsilon$.
In a similar way,
$\mu(g^{-1}\widetilde Y)=\mu_1^{k_\epsilon}(g_{k_\epsilon}^{-1}Y_{k_\epsilon}^c)<\epsilon$.
It follows that
$$
 \int_{\widetilde Y}(f\circ g_n^{-1}\cdot f)(x_1,\dots,x_n)\sqrt{\frac{d\mu\circ g^{-1}}{d\mu}}(x)\,d\mu(x)< C^2\epsilon.
$$
This yields that
 $$
 \aligned
\Bigg| \langle\kappa_\mu(g) (f\otimes 1),f\otimes 1\rangle
 &-
 \int_{Y_{k_\epsilon}}(f\circ g_n^{-1}\cdot f)(y_1,\dots,y_n)\sqrt{\frac{d(\mu^{k_\epsilon}_1\circ g_{k_\epsilon}^{-1})}{d\mu_1^{k_\epsilon}}}(y)\,d\mu_1^{k_\epsilon}(y)\Bigg| \\
 &\le C^2\epsilon.
 \endaligned
% \tag7-1
 $$
 where $(y_1,\dots, y_n)$ are the $n$ first coordinates of $y\in X_1^{k_\epsilon}$.
 Let $\Phi:=\sqrt{\frac{d\mu_1^{k_\epsilon}}{d\nu_1^{k_\epsilon}}}$.
 Then $\Phi$ is a strictly positive function defined on $X_1^{k_\epsilon}$.
\comment

We now let, for each $y=(y_1,\dots,y_{k_\epsilon})\in X_1^{k_\epsilon}$,
 $$
 \widetilde f(y):=f(y_1,\dots,y_n)\cdot \Phi(y)\cdot 1_{Y_{k_\epsilon}}(y).
 $$
 Then $ \widetilde f$ is a mapping from $X_1^{k_\epsilon}$ to $\Bbb C$ and
 $$
 \multline
 \int_{Y_{k_\epsilon}}f\circ g_n^{-1}(y_1,\dots,y_n)f(y_1,\dots,y_n)\sqrt{\frac{d(\mu^{k_\epsilon}_1\circ g_{k_\epsilon}^{-1})}{d\mu_1^{k_\epsilon}}}(y)\,d\mu_1^{k_\epsilon}(y)\\
 =
 \int_{Y_{k_\epsilon}}\widetilde f\circ g_n^{-1}(x_1,\dots,x_n)\widetilde f(x_1,\dots,x_n)\sqrt{\frac{d(\nu^{k_\epsilon}_1\circ g_{k_\epsilon}^{-1})}{d\nu_1^{k_\epsilon}}}(y)\,d\nu_1^{k_\epsilon}(y)\\
  =\langle\kappa_\nu(g) (\widetilde f\otimes 1),\widetilde f\otimes 1\rangle
 \endmultline
 $$

 \endcomment
 We define a mapping $ \widetilde f:X\to\Bbb C$ by setting
  $$
 \widetilde f(x):=f(x_1,\dots,x_n)\cdot \Phi(x_1,\dots,x_{k_\epsilon})\cdot 1_{Y_{k_\epsilon}}(x).
 $$
 Then
 $$
 \align
  \langle\kappa_\nu(g) (\widetilde f\otimes 1),\widetilde f\otimes 1\rangle
& =
 \int_{Y_{k_\epsilon}}
( \widetilde f\circ g_{k_\epsilon}^{-1}\cdot\widetilde f)(y)
\sqrt{\frac{d(\nu^{k_\epsilon}_1\circ g_{k_\epsilon}^{-1})}{d\nu_1^{k_\epsilon}}}(y)\,d\nu_1^{k_\epsilon}(y)\\
&=
 \int_{Y_{k_\epsilon}\cap g_{k_\epsilon}Y_{k_\epsilon}}
( f\circ g_n^{-1}\cdot f)(y_1,\dots,y_n)
\sqrt{\frac{d(\mu^{k_\epsilon}_1\circ g_{k_\epsilon}^{-1})}{d\mu_1^{k_\epsilon}}}(y)\,d\mu_1^{k_\epsilon}(y)\\
&=\int_{Y_{k_\epsilon}}-\int_{Y_{k_\epsilon}\setminus g_{k_\epsilon}Y_{k_\epsilon}}.
\endalign
 $$
By the Cauchy-Schwarz inequality,
 $$
 \bigg|\int_{Y_{k_\epsilon}\setminus g_{k_\epsilon}Y_{k_\epsilon}}\bigg|\le C^2
 \sqrt{\int_{Y_{k_\epsilon}\setminus g_{k_\epsilon}Y_{k_\epsilon}}d\mu_1^{k_\epsilon}\,
 \int_{g_{k_\epsilon}^{-1}Y_{k_\epsilon}\setminus Y_{k_\epsilon}}d\mu_1^{k_\epsilon}}
 \le C^2\epsilon.
 $$
 Therefore,
 $$
\max_{g\in F} |\langle\kappa_\mu(g) (f\otimes 1),f\otimes 1\rangle-
 \langle\kappa_\nu(g) (\widetilde f\otimes 1),\widetilde f\otimes 1\rangle
 |<2C^2\epsilon.
 $$
 Since the linear subspace $\bigcup_{n=1}^\infty\{f\otimes 1\mid f:X_1^n\to\Bbb C\}$
 is dense in $L^2(X,\mu)$, it follows that $\kappa_\mu\prec\kappa_\nu$.
% In a similar way, $\kappa_\nu\prec\kappa_\mu$.
% Hence, $\kappa_\mu$ and $\kappa_\nu$ are weakly equivalent.
 \qed
\enddemo

Part (v) of Main Result follows from Theorem~7.1.

\proclaim{Corollary 7.2} Under the condition of Theorem~6.5, the Koopman representations $\kappa_{\mu_\omega}$, $\omega\in\{0,1\}^\Bbb N$, are pairwise weakly equivalent.
\endproclaim

\comment

\head 4. What if $G$ is weakly branch but not saturated
\endhead

If $G$ is not saturated but weakly branch then consider $G_O$ where $O$ is a cylinder.
Then $(O,\mu\restriction O, G_O)$ is not necessarily ergodic.
However, it is easy to describe an ergodic decomposition of this system.
Let $G_O^{(n)}\subset \text{Bij}\,X_1^n$ is the projection of $G_O$ to $X_1^n$.
There is a partition $P_1$ into $G_O^{(1)}$-orbits.
Then for each $p_1\in P_1$, the set $p_1\times X_2$ is $G_O^{(2)}$-invariant

let $P_{2,p_1}$ be the partition of $X_2$ into

\endcomment

\enddocument